\documentclass[11pt]{amsart}
\usepackage{amsmath}
\usepackage{amssymb}
\usepackage{amscd}
\usepackage[latin2]{inputenc}
\usepackage[T1]{fontenc}
\usepackage{pstricks}
\usepackage{pstricks-add}
\usepackage{pst-node}
\usepackage{pst-coil}
\usepackage{graphicx}

\newcommand{\CC}{\mathbb{C}}
\newcommand{\intfrac}[2]{\genfrac{\lfloor}{\rfloor}{}{1}{#1}{#2}}

\newcommand{\scalepic}{1.0}

\newcommand{\zhand}{$0$--handle{}}
\newcommand{\ohand}{$1$--handle{}}
\newcommand{\ohands}{$1$--handles{}}
\newcommand{\thand}{$2$--handle{}}
\newcommand{\thands}{$2$--handles{}}

\newcommand{\nz}{n^*_\zeta}
\newcommand{\nnsz}{n_\zeta}
\newcommand{\sz}{\sigma^*_\zeta}
\newcommand{\snsz}{\sigma_\zeta}
\newcommand{\LS}{L^{sing}}
\newcommand{\pg}{{p_g}}

\newcommand{\schempic}{\begin{pspicture}(-8,-0.5)(8,2.5)
\psscalebox{0.80}{
\pscircle[fillstyle=solid,fillcolor=lightgray,opacity=0.9](-6,2){0.75}%
\rput(-6,3){\psscalebox{\scalepic}{$S_{-}$}}
\pscircle[fillstyle=solid, fillcolor=black](-5,2){0.03}%
\rput(-5,2.2){\psscalebox{\scalepic}{$z_k$}}
\psline{->}(-4.5,2)(-3,2)%
\rput(-3.75,2.2){\psscalebox{\scalepic}{Step 1}}
\pscircle[fillstyle=solid,fillcolor=lightgray,opacity=0.9](-2,2){0.75}%
\rput(-2,3){\psscalebox{\scalepic}{$S_{-}$}}
\pscircle[fillstyle=solid,fillcolor=lightgray,opacity=0.9](-1,2){0.27}
\pscircle[fillstyle=solid, fillcolor=black](-1,2){0.03}%
\psline{->}(-0.25,2)(1,2)\rput(0.4,2.2){\psscalebox{\scalepic}{Step 2}}%
\rput(2,3){\psscalebox{\scalepic}{$S_{-}$}}
\pscircle[fillstyle=solid, fillcolor=lightgray,opacity=0.9, linestyle=dotted](3,2){0.5}
\pscircle[fillstyle=solid, fillcolor=lightgray,opacity=0.9, linestyle=dotted](2,2){0.75}
\psarc(2,2){0.75}{30}{330}
\psarc(3,2){0.5}{222}{138}
\pscircle[fillstyle=solid, fillcolor=black](3,2){0.03}%
\psline{->}(-6.7,0)(-5.2,0)\rput(-6.1,0.2){\psscalebox{\scalepic}{Step 3}}
\rput(-5.3,-1){\psscalebox{\scalepic}{$S(0,\tilde{r})$}}
\pscircle[fillstyle=solid, fillcolor=lightgray,opacity=0.9, linestyle=dotted](-3,0){0.5}
\pscircle[fillstyle=solid, fillcolor=lightgray,opacity=0.9, linestyle=dotted](-4,0){1.12}
\psarc(-4,0){1.12}{28}{332}
\psarc(-3,0){0.5}{266}{94}
\rput(-3,0.2){\psscalebox{0.99}{$z_k$}}
\pscircle[fillstyle=solid,fillcolor=black](-3,0){0.03}
\psline{->}(-2,0)(-0.5,0)\rput(-1.3,0.2){\psscalebox{\scalepic}{Step 4}}
\pscircle[fillstyle=solid,fillcolor=lightgray,opacity=0.9](1,0){1.12}
\rput(0.5,-1.25){\psscalebox{\scalepic}{$S_+$}}
\pscircle[fillstyle=solid,fillcolor=black](2.0,0){0.03}
}
\end{pspicture}}

\newcommand{\nextpictureonC}{\begin{pspicture}(-5,-1.5)(5,1.5)
\begin{psclip}{\pscircle[linestyle=none](-3.0,0.0){1.2}}
\psellipse[fillcolor=lightgray,fillstyle=solid,linewidth=1.2pt](-3.0,1.5)(0.8,1.3)
\pscircle[fillcolor=lightgray,fillstyle=solid,linewidth=1.2pt](-3.0,-1.5){1.3}
\end{psclip}%
\begin{psclip}{\pscircle[linestyle=none](2.0,0.0){1.2}}
\psellipse[fillcolor=lightgray,fillstyle=solid,linestyle=solid,linecolor=black,linewidth=1.2pt](2.0,1.5)(0.8,1.3)
\pscircle[fillcolor=lightgray,fillstyle=solid,linecolor=black,linewidth=1.2pt](2.0,-1.5){1.9}
\end{psclip}%
\begin{psclip}{\psellipse[linestyle=none](2.0,1.5)(0.8,1.3)}
\pscircle[linecolor=lightgray,linewidth=2.0pt](2.0,-1.5){1.9}
\end{psclip} %
\begin{psclip}{\pscircle[linestyle=none](2.0,-1.5){1.9}}
\psellipse[linecolor=lightgray,linewidth=2.0pt](2.0,1.5)(0.8,1.3)
\end{psclip}%
\pscircle(1.9,0.0){0.9}
\pscircle(-3.1,0.0){0.9}
\rput(-1.85,0.0){\psscalebox{0.8}{$D_{ja}$}}
\rput(3.0,0.4){\psscalebox{0.8}{$D_{ja}$}}
\rput(-3.1,1.4){\psscalebox{0.8}{$\Gamma$}}
\rput(1.9,1.4){\psscalebox{0.8}{$\Gamma$}}
\rput(-3.1,-1.4){\psscalebox{0.8}{$R_s$}}
\rput(1.9,-1.4){\psscalebox{0.8}{$R_s$}}
\end{pspicture}}

\newcommand{\toymodel}{
\begin{pspicture}(-6,-1.5)(6,1.5)
\pscircle[fillcolor=lightgray, fillstyle=solid, opacity=0.9,linestyle=none](-4,0){1.0}
\pscircle[fillcolor=lightgray, fillstyle=solid, opacity=0.9,linestyle=none](-2.1,0){1.0}
\psellipse[linestyle=dotted](-4,0)(1.0,0.2)
\psellipse[linestyle=dotted](-2.1,0)(1.0,0.2)
\psarc[linecolor=black](-4,0){1.0}{20}{340}
\psarc[linecolor=black](-2.1,0){1.0}{200}{160}
\pspolygon[fillstyle=none, opacity=0.7, linewidth=0.3pt, linecolor=darkgray](-5.5,0.5)(-1,0.5)(-0.5,0.8)(-5,0.8)
\psellipse[linecolor=black](-4,0.65)(0.75,0.15)
\psellipse[linecolor=black](-2.1,0.65)(0.75,0.15)
\rput(-4.0,-1.2){\psscalebox{0.8}{$B_1$}}
\rput(-2.1,-1.2){\psscalebox{0.8}{$B_2$}}
\rput(-5.25,0.9){\psscalebox{0.8}{$C$}}
\pscircle[fillcolor=lightgray, fillstyle=solid, opacity=0.9,linestyle=none](2,0){1.0}
\pscircle[fillcolor=lightgray, fillstyle=solid, opacity=0.9,linestyle=none](3.1,0){1.0}
\psellipse[linestyle=dotted](2,0)(1.0,0.2)
\psellipse[linestyle=dotted](3.1,0)(1.0,0.2)
\pspolygon[fillstyle=none, opacity=0.7, linewidth=0.3pt, linecolor=darkgray](0.5,0.5)(4.5,0.5)(5,0.8)(1,0.8)
\psarc[linecolor=black](2,0){1.0}{57}{303}
\psarc[linecolor=black](3.1,0){1.0}{237}{123}
\psellipticarc[linecolor=black](2,0.65)(0.75,0.15){40}{320}
\psellipticarc[linecolor=black](3.1,0.65)(0.75,0.15){220}{140}
\rput(2.0,-1.2){\psscalebox{0.8}{$B_1$}}
\rput(3.1,-1.2){\psscalebox{0.8}{$B_2$}}
\rput(0.75,0.9){\psscalebox{0.8}{$C$}}
\end{pspicture}
}

\newcommand{\gluepic}{
\begin{pspicture}(-6,-1.2)(6,2.2)
\rput(0,0){\psscalebox{0.6}{%
\psellipticarc(-5,-3.5)(0.8,3){30}{150}
\psellipticarc[linewidth=3pt](-5,-3.5)(0.8,3){70}{110}
\psellipticarc(-5,4.7)(0.8,3){215}{325}
\psellipticarc[linewidth=3pt](-5,4.7)(0.8,3){250}{290}
\psarc(-5,1){0.5}{315}{295}
\psarc[linewidth=2pt](-5,1){0.5}{70}{110}
\psarc(-5,0.2){0.5}{135}{115}
\psarc[linewidth=2pt](-5,0.2){0.5}{250}{290}
\psline[linewidth=2pt]{->}(-4,0.6)(-2.3,0.6)
\psellipticarc(-1.5,-4.5)(0.7,5.5){27}{106}
\psellipticarc(-1.5,-4.5)(0.7,5.5){116}{153}
\psellipticarc(-1.5,5.7)(0.7,5.5){209}{285}
\psellipticarc(-1.5,5.7)(0.7,5.5){295}{331}
\psellipticarc(2,-3.8)(0.6,3){37}{143}
\psellipticarc[linewidth=3pt](2,-3.8)(0.6,3){70}{110}
\psellipticarc(2,5)(0.6,3){220}{320}
\psellipticarc[linewidth=3pt](2,5)(0.6,3){254}{286}
\psarc(2,1){0.6}{330}{250}
\psarc[linewidth=3pt](2,1){0.6}{75}{105}
\psarc(2,0.2){0.6}{150}{70}
\psarc[linewidth=3pt](2,0.2){0.6}{255}{285}
\psbezier(1.62,0.62)(1.82,0.82)(1.92,0.72)(1.82,0.72)
\psbezier(1.82,0.75)(2.02,0.72)(2.02,0.52)(2.22,0.52)
\psbezier(2.22,0.52)(2.27,0.52)(2.32,0.52)(2.37,0.57)
\psbezier(1.8,0.44)(1.95,0.44)(1.85,0.4)(1.97,0.57)
\psline[linewidth=2pt]{->}(2.8,0.6)(4.5,0.6)
\psbezier(5.12,0.62)(5.32,0.82)(5.32,0.72)(5.35,0.72)
\psbezier(5.32,0.75)(5.52,0.72)(5.52,0.52)(5.72,0.52)
\psbezier(5.72,0.52)(5.77,0.52)(5.82,0.52)(5.87,0.57)
\psellipticarc(5.5,-4.3)(1.1,5.5){25}{80}
\psellipticarc(5.5,-4.3)(1.1,5.5){120}{155}
\psellipticarc(5.5,5.6)(1.1,5.5){208}{260}
\psellipticarc(5.5,5.6)(1.1,5.5){300}{332}
\psbezier(5.7,1.1)(5.6,1.2)(5.5,1.0)(5.5,0.8)
\psbezier(5.3,0.2)(5.4,0.1)(5.5,0.45)(5.5,0.55)}}
\end{pspicture}
}

\newcommand{\whykappa}{\begin{pspicture}(-6,-2)(6,2)
\rput(0,0){\psscalebox{0.8}{
\psarc(-4.5,0){1.5}{17}{343}\psarc(-2.85,0){0.5}{242}{118}
\pscircle[fillcolor=black,fillstyle=solid](-4.5,0){0.05}
\pscircle[fillcolor=black,fillstyle=solid](-2.85,0){0.05}
\psellipse[linestyle=dotted](-3.05,0)(0.1,0.43)
\psline[linecolor=gray](-3.90,1.5)(-1.80,-1.5)\rput(-1.50,-1.5){\psscalebox{0.80}{$G_3$}}
\psline[linecolor=gray](-4.20,1.0)(-1.50,-1.0)\rput(-1.20,-1.0){\psscalebox{0.80}{$G_2$}}
\psline[linecolor=gray](-4.20,-1.0)(-1.50,1.0)\rput(-1.20,1.0){\psscalebox{0.80}{$G_1$}}
\rput(-4.5,-1.7){\psscalebox{0.8}{$S_-$}}
\rput(-2.85,-0.7){\psscalebox{0.8}{$S_\mu$}}
\rput(-3.8,1.8){Good:}
\psarc(1.5,0){1.5}{10}{350}\psarc(3.15,0){0.3}{240}{120}
\pscircle[fillcolor=black,fillstyle=solid](1.5,0){0.05}
\pscircle[fillcolor=black,fillstyle=solid](3.15,0){0.05}
\psellipse[linestyle=dotted,dotsep=1pt](2.97,0)(0.08,0.23)
\psline[linecolor=gray](3.60,-1.5)(2.70,1.5)\rput(3.90,-1.5){\psscalebox{0.80}{$G_3$}}
\psline[linecolor=gray](4.80,1.0)(1.50,-1.0)\rput(5.10,-1.0){\psscalebox{0.80}{$G_2$}}
\psline[linecolor=gray](4.80,-1.0)(1.50,1.0)\rput(5.10,1.0){\psscalebox{0.80}{$G_1$}}
\rput(1.5,-1.7){\psscalebox{0.8}{$S_-$}}
\rput(3.15,-0.5){\psscalebox{0.8}{$S_\mu$}}
\rput(3.2,1.8){Bad:}}}
\end{pspicture}}
\newcommand{\notationonAj}{\begin{pspicture}(-5,-3)(5,3)
\begin{psclip}{\pscircle(0,0){2.5}}
\rput{0}(0,0){\psellipse[fillcolor=black,fillstyle=crosshatch*,hatchcolor=white,hatchwidth=0.8pt, hatchsep=0.5pt,opacity=1.0](3,0)(1.8,0.5)
\pscircle[fillcolor=black,fillstyle=solid](1.2,0.0){0.05}}
\rput{120}(0,0){\psellipse[fillcolor=black,fillstyle=crosshatch*,hatchcolor=white,hatchwidth=0.8pt, hatchsep=0.5pt](3,0)(1.8,0.5)
\pscircle[fillcolor=black,fillstyle=solid](1.2,0.0){0.05}}
\rput{240}(0,0){\psellipse[fillcolor=black,fillstyle=crosshatch*,hatchcolor=white,hatchwidth=0.8pt, hatchsep=0.5pt](3,0)(1.8,0.5)
\pscircle[fillcolor=black,fillstyle=solid](1.2,0.0){0.05}
\pscircle[fillcolor=lightgray,opacity=0.4,fillstyle=vlines](1.2,0.0){0.3}}
\rput(0.9,0.0){\psscalebox{0.8}{$y_{s1}$}}
\psarc(0,0){1.2}{90}{150}
\rput(0,1.35){\psscalebox{0.8}{$\partial R_{s2}$}}
\rput(0.0,-1.1){\psscalebox{0.8}{$D_{s3}$}}
\pscircle[fillstyle=solid,fillcolor=black](0,0){0.05}
\end{psclip}
\rput(2.7,0){\psscalebox{0.8}{$\Gamma$}}
\rput(-1.4,2.3){\psscalebox{0.8}{$\Gamma$}}
\rput(-1.4,-2.3){\psscalebox{0.8}{$\Gamma$}}
\pscircle[fillcolor=darkgray,fillstyle=solid,opacity=0.4](0,0){0.5}\rput(0.0,0.7){\psscalebox{0.8}{$R_{\mu}$}}
\rput(0.0,2.7){\psscalebox{0.8}{$R_\lambda$}}
\end{pspicture}}
\newcommand{\stepfour}{\begin{pspicture}(-5,-2.5)(5,2.5)%
\pscircle[fillcolor=black,fillstyle=solid](0.0,0){0.05}\rput(0,0.2){\psscalebox{0.8}{$z_k$}}%
\psarc[linewidth=2.0pt](0,0){1.5}{270}{90}\rput(2.1,0.0){\psscalebox{0.8}{$\phi_0^4(W_S)$}}%
\pscircle(0,0){1.5}\rput(-2.2,0.0){\psscalebox{0.8}{$S(z_k,\lambda\varepsilon)$}}%
\begin{psclip}{\pscircle[linestyle=none](0.0,0.0){2.3}}%
\pscircle[linewidth=0.5pt](-5,0){5.22}%
\end{psclip}%
\begin{psclip}{\pscircle[linestyle=none](0.0,0.0){1.5}}%
\pscircle[linewidth=2.0pt](-5,0){5.22}%
\pscircle[linestyle=dashed](-2,0){2.5}%
\pscircle[linestyle=dashed](-0.8,0){1.7}%
\pscircle[linestyle=dashed](-0.4,0){1.55}%
\end{psclip}%
\psline(-2,-2)(2,2)\rput(1.8,2.1){\psscalebox{0.8}{$Z_1$}}
\psline(3,-2)(-3,2)\rput(2.8,-1.6){\psscalebox{0.8}{$Z_2$}}
\rput(0.35,2.0){\psscalebox{0.8}{$S(0,\tilde{r})$}}
\end{pspicture}}

\DeclareMathOperator{\re}{Re}

\newtheoremstyle{captionstyle}{6pt}{2pt}{\upshape}{0em}{\scshape}{.}{ }{}%

\theoremstyle{definition}
\newtheorem{example}{Example}[section]
\newtheorem{definition}[example]{Definition}
\theoremstyle{captionstyle}
\newtheorem{psfiguretheo}{Figure}

\theoremstyle{plain}
\newtheorem{lemma}[example]{Lemma}
\newtheorem{proposition}[example]{Proposition}
\newtheorem{theorem}[example]{Theorem}

\newtheorem{corollary}[example]{Corollary}

\theoremstyle{remark}
\newtheorem*{acknowledgements}{Acknowledgements}
\newtheorem{remark}[example]{Remark}

\newtheorem{convention}[example]{Convention}
\numberwithin{equation}{section}

\newenvironment{psfigure}%
{\footnotesize
\begin{list}{}%
{\setlength{\leftmargin}{1cm}\setlength{\rightmargin}{1cm}}%
\item[]\begin{psfiguretheo}
}%
{\end{psfiguretheo}\end{list}}%
\begin{document}
\title[Morse theory]{Morse theory for plane algebraic curves}
\author{Maciej Borodzik}
\address{Institute of Mathematics, University of Warsaw, ul. Banacha 2,
02-097 Warsaw, Poland}
\email{mcboro@mimuw.edu.pl}
\date{April 04, 2011}
\thanks{Supported by Polish KBN Grant No 2 P03A 010 22 and Foundation for Polish Science}
\begin{abstract} 
We use Morse theorical arguments to study algebraic curves in $\mathbb{C}^2$. We take an algebraic 
curve $C\subset\mathbb{C}^2$ and intersect it with
spheres with fixed origin and growing radii. We explain in detail how the embedded type of the intersection
changes if we cross a singular point of $C$. Then we apply link invariants such as Murasugi's signature and Tristram--Levine
signature to obtain informations about possible singularities of the curve $C$ in terms of its topology. 
\end{abstract}
\maketitle

\section{Introduction}
By a plane algebraic curve we understand a set 
\[C=\{(w_1,w_2)\in\mathbb{C}^2\colon F(w_1,w_2)=0\},\] 
where $F$ is an irreducible polynomial.
Let $\xi=(\xi_1,\xi_2)\in\mathbb{C}^2$, and $r\in\mathbb{R}$ be positive. If the intersection of $C$ with a 3-sphere 
$S(\xi,r)$
is transverse, it is a link in $S(\xi,r)\simeq S^3$. We denote it by $L_r$.

If $\xi$ happens to be a singular point of $C$ and $r$ is sufficiently small, $L_r$ is a link of a plane curve singularity 
of $C$ at $\xi$. 
On the other hand, for any $\xi\in\mathbb{C}^2$ and
for any  sufficiently large $r$, $L_r$ is the link of $C$ at infinity.

Links of plane curve singularities have been perfectly understood for almost thirty years (see \cite{EN} for topological 
or \cite{Wall} for algebro-geometrical approach). 
Possible links at infinity are also well
described (see \cite{Neu2, NeRu}). The most difficult case to study, 
as it was pointed out in a beautiful survey \cite{Rud0}, is the intermediate
step, i.e. possible links $L_r$ for $r$ neither very small nor very large.

Our idea is to study the differences between the links of singularities of a curve and its link at infinity via Morse theory:
we begin with $r$ small and let it grow to infinity. The isotopy type of the link changes, when we pass through critical points.
If $C$ is smooth, the theory is classical (see e.g. \cite[Chapter V]{Ka} or \cite{Mil}), yet if $C$ has singular points,
the analysis requires more care and is a new element in the theory.

To obtain numerical relations we apply some knot invariants. Namely, we study changes of Murasugi's signature in detail and then pass
to Levine--Tristram signatures, which give a new set of information. Our choice is dictated by the fact, that these invariants are well behaved under
the one handle addition (this is Murasugi's Lemma, see Lemma~\ref{mur}). From a knot theoretical point of view, Morse theory provides
inequalities between signatures, which are very closely related to those in \cite{Kaw,Kaw2} (cf. Corollary~\ref{betti} and a discussion below it).
What is important, are the applications in algebraic geometry. In this paper we show only a few of them. 
First of all, we present an elementary proof of Corollary~\ref{cormain}. The only known proof up to now \cite{BZ2,BZ3} relies heavily
on algebraic geometry techniques. This result is of interest not only for algebraic geometers, but also in the theory of bifurcations of ODE's (see
\cite{ChL,BZ3} and references therein). We also reprove Varchenko's estimate on the number of cusps of a degree $d$ curve in $\mathbb{C}P^2$ (see
Corollary~\ref{cor:2372}). Corollary~\ref{ratknot} and Lemma~\ref{ratknot2} show also a different, completely new application of our method.
We refer to \cite{Bo} for a brand new application in studying deformations of singularities of plane curves. 

We want also to point out that the methods developped in this article have been used in \cite{BN2} to
show various semicontinuity results for singularities of plane curves --- including establishing a relationship between spectrum of a polynomial 
in two variables at infinity and spectra of singular points of one of its fibers --- in a purely topological way. The application
of (generalized) Tristram--Levine signatures in higher dimensional singularity theory is also possible, 
even though the details somehow differ from those developped in the present paper. This latter work is in progress.

Albeit Tristram--Levine signatures turn out to be an important tool of extracting data about plane curves, it is surely not the only one.
One of the main messages of the article is that \emph{any knot cobordism invariant can be used to obtain global informations about possible
singularities which may occur on a plane curves}. Altough the $s$ invariant of Rasmussen \cite{Ras} and the $\tau$ invariant of Ozsv\'ath--Szabo
\cite{OS} apparently do not give any new obstructions (they are equal to the four genus for positive knots) and Peters' invariant \cite{Pts} seem to be
very much related to the Tristram--Levine signature at least for torus knots, but the author is convinced that the 
application of full Khovanow homology in this context will lead to brand new discoveries in the theory of plane curves.
  
\begin{convention}
Throughout the paper we use standard Euclidean, metric on $\mathbb{C}^2$. $B(\xi,r)$ denotes the ball
with centre $\xi$ and radius $r$. We may assume, to be precise, that it is a closed ball, but we never appeal to this fact.
The boundary of the ball $B(\xi,r)$ is the sphere denoted $S(\xi,r)$.
\end{convention}

\section{Handles related to singular points}
Let $C$ be a plane algebraic curve given by equation $F=0$, where $F$ is a reduced polynomial. Let $\xi\in\mathbb{C}^2$.
Let $z_1,\dots,z_n$ be all the points of $C$ such that either $C$ is not transverse to $S(\xi,||z_k-\xi||)$
at $z_k$, or $z_k$ is a singular point of $C$. We shall call them critical points. Let
\[\rho_k=||z_k-\xi||.\]
We order $z_1,\dots,z_n$ in such a way that $\rho_1\le \rho_2\le \dots \le \rho_n$. 
We shall call $\rho_k$'s \emph{critical values}.
We shall pick a \emph{generic} $\xi$ which
means that
\begin{itemize}
\item[(G1)] $\rho_1<\rho_2<\dots<\rho_n$, i.e. at each level set of the distance function 
\begin{equation}\label{eq:dist}
g=g_\xi(w_1,w_2)=|w_1-\xi_1|^2+|w_2-\xi_2|^2
\end{equation}
restricted to $C$ there is at most one critical point (this is not a very serious restriction and it is put here rather for convenience).
\item[(G2)] If $z_k$ is a smooth point of $C$, then $g|_C$ is of Morse type near $z_k$.
\item[(G3)] If $z_k$ is a singular point of $C$, we assume the condition \eqref{eq:tan2} holds.
\end{itemize}
Generic points always exist. Obviously G3 and G1 are open-dense conditions. For G2 see e.g. \cite[Theorem~6.6]{Mil}.

We want to point out that we assume here tacitly, that the overall number of critical points is finite. This follows from the 
algebraicity of the curve $C$ (see Remark~\ref{rem:fin}).
if $C$ is not algebraic, this does not hold automatically, because even the number of singular points of $C$ can be infinite and the
link at infinity hard to define at all,  consider e.g.
a curve $\{(z_1,z_2)\in\mathbb{C}^2\colon z_1\sin z_2=0\}$. Using methods of \cite[Proposition 2]{FGR} one can produce
other amusing, albeit not explicit, examples.

\begin{remark}\label{rem:smooth1}
From the condition G3 we see in particular that if $\xi$ does not lie on $C$, then $z_1$ is a smooth point of $C$. Indeed, $g|_C$ attains
local minimum of $z_1$, so the tangent space $T_{z_1}C$ is not transverse to $T_{z_1}S(\xi,\rho_1)$. If $z_1$ is not smooth, this violates G3.
\end{remark}

It is well known that, if $r_1$ and $r_2$ are in the same interval $(\rho_k,\rho_{k+1})$ then links
$L_{r_1}$ and $L_{r_2}$ are isotopic, where
\[L_r=C\cap S(\xi,r)\subset S(\xi,r).\]
Next definition provides very handy language.
\begin{definition}
Let $\rho_k$ be a critical value. The links $L_{\rho_k+}$ and $L_{\rho_k-}$ (or, if there is no risk of confusion, 
just $L_+$, $L_-$) are
the links $L_{\rho_k+\varepsilon}$ and $L_{\rho_k-\varepsilon}$
with $\varepsilon>0$ such that $\rho_k+\varepsilon<\rho_{k+1}$ and
$\rho_k-\varepsilon>\rho_{k-1}$. We shall say, informally, that the change from $L_-$ to $L_+$
is a crossing or a passing through a singular point $z_k$.
\end{definition}
\begin{lemma}
Assume that $z_k$ is a smooth point of $C$. Then $L_{\rho_k+}$ arises
from $L_{\rho_k-}$ by addition of a \zhand{}, an \ohand{} or a \thand{} according to the Morse index
at $z_k$ of the distance function $g$
restricted to $C$.
\end{lemma}
A \zhand\ corresponds to adding an unlinked unknot to the link. A \thand\
corresponds to deleting an unlinked unknot. The addition of a \ohand\ is a hyperbolic operation, which we now define.
\begin{definition}[see \expandafter{\cite[Definition~12.3.3]{Kaw-book}}]\label{hyper}
Let $L$ be a link with components $K_1,\dots,K_{n-1},K_n$. 
Let us join the knots $K_{n-1}$ and $K_n$ by a band, so as to obtain a knot $K'$.
Let $L'=K_1\cup\dots\cup K_{n-2}\cup K'$. We shall then say, that $L'$ is obtained from $L$ by
a \emph{hyperbolic transformation}.
\end{definition}
The hyperbolic transformation depends heavily on the position of the band, for example, by adding a band to a Hopf link we can obtain
a trivial knot, but also a trefoil and, in fact, infinitely many different knots.

\begin{remark}\label{rem:smooth2}
Assume again that $\xi\not\in C$. We know that $z_1$ is a smooth point. As for $r<\rho_1$ the link $L_r$ is empty and for $r>r_1$ it is not, the
first handle must be a birth. In particular, for $r\in(\rho_1,\rho_2)$ the link $L_r$ is an unknot.
\end{remark}

\begin{lemma}
If $C$ is a complex curve, there are no \thands.
\end{lemma}
\begin{proof}
A \thand\ corresponds to a local maximum of a distance function \eqref{eq:dist}
restricted to $C$. The functions $w_1-\xi_1$ and $w_2-\xi_2$ are holomorphic on $C$, hence 
$|w_1-\xi_1|^2+|w_2-\xi_2|^2$ is subharmonic on $C$, and as such, it does not have any local maxima on $C$.
\end{proof}
\ohands\ might occur in three forms. 

\begin{definition}
Let $C_-=C\cap B(\xi,\rho_k-\varepsilon)$.
A \ohand{} attached to two different connected components of the normalization of $C_-$
is called a \emph{join}. A \ohand{} attached to a single component of the normalization of $C_-$ but to two different
components of $L_-$, is called a \emph{marriage}. And finally, if it is attached to a single component
of $L_-$, it is called a \emph{divorce}.
\end{definition}

If the point $z_k$ is not smooth, the situation is more complicated.
\begin{definition}
The \emph{multiplicity} of a singular point $z$ of $C$ is the local intersection index
of $C$ at $z$ with a generic line passing through $z$.
\end{definition}
\begin{proposition}\label{singhand}
Let $z_k$ be a singular point of $C$ with multiplicity $p$.
Let $\LS$ 
be the link of the singularity at $z_k$. 
Then $L_+$ $(=L_{\rho_k+})$ can be obtained from the disconnected sum of $L_-$ $(=L_{\rho_k-})$ with $\LS$
by adding $p$ \ohands. 
\end{proposition}
\begin{proof}
This is the most technical and difficult proof in the article. First we shall introduce the notation, then we shall outline the proof, which in turn 
will consist in four steps.  

\underline{\emph{Introducing the notation.}} Up to an isometric coordinate change we can 
assume that $\xi=(0,0)$ and $z_k=(\rho_k,0)$. 

Let $G_1,\dots,G_b$ be the branches of $C$ at $z_k$.
By Puiseux theorem (see e.g. \cite[Section 2]{Wall}), each branch $G_j$ can be locally parametrized in a Puiseux expansion
\begin{equation}\label{eq:param}
w_1=\rho_k-\beta_j\tau^{p_j},\,\,w_2=\alpha_j\tau^{p_j}+\dots,\,\,\,\,\,\, \tau\in\mathbb{C},\,\,|\tau|\ll 1,
\end{equation}
i.e. it is a topological disk. Let $\psi_j\colon\{|\tau|\ll 1\}\to\mathbb{C}^2$ be the parametrization given by \eqref{eq:param}.

The (generalised) tangent line
to $G_j$ at $z_k$ is the line $Z_j$ defined by
\begin{equation}\label{eq:tangentlines}
Z_j=\{(w_1,w_2)\in \mathbb{C}^2\colon \alpha_j(w_1-\rho_k)+\beta_j w_2=0.\}
\end{equation}
The tangent space to $C$ at $z_k$ is then the union of lines $Z_1,\dots,Z_b$.
By genericity of $\xi$ we may assume that
\begin{equation}\label{eq:tan2}
\alpha_j\beta_j\neq 0\text{ for any $j$}. 
\end{equation}
This means that nether the line $\{(w_1,w_2)\colon w_1-\rho_k=0\}$ nor $\{w_2=0\}$ is tangent to $C$ at $z_k$. In other words, we can choose $\varepsilon$, 
$\lambda$ and $\mu$ in such way, that the following conditions are satisfied.
\begin{itemize}
\item[(S1)] The intersection of each tangent line $Z_j$ with $S(0,\rho_k-\varepsilon)$ is non-empty (we use $\beta_j\neq 0$);
\item[(S2)] The intersection $B(0,\rho_k-\varepsilon)\cap B(z_k,\mu\varepsilon)$ is non-empty and \emph{omits} each tangent line $Z_j$ 
(i.e. $\mu>1$, $\mu$ is very close to $1$ and we use $\alpha_j\neq 0$);
\item[(S3)] The two-sphere $S(0,\rho_k-\varepsilon)\cap S(z_k,\lambda\varepsilon)$ is not disjoint with $Z_j$ (this is a refinement of (S1));
\item[(S4)] $\lambda\varepsilon$ is sufficiently small (in the sense which will be made precise later);
\item[(S5)] In particular, if we choose
\[\tilde{r}=\sqrt{\rho_k^2+\lambda^2\varepsilon^2},\]
then $z_k$ is the only point at which the intersection of $C$ with $S(0,r)$ is not transverse, for $r\in[\rho_k-\varepsilon,\tilde{r}]$.
\end{itemize}

It is important to show that the two conditions $\alpha_j\neq 0$ and $\beta_j\neq 0$ are of different nature. Namely, if for some $j$, $\beta_j=0$,
the proposition fails. On the other hand, the condition $\alpha_j\neq 0$ is used only to make the exposition clearer and easier to understand. The
proof given below works if for some $j$, $\alpha_j=0$, but we would have use less transparent arguments in two places.

Let us define the following sets:
\begin{align*}
B_{-}&=B(0,\rho_k-\varepsilon)&B_+&=B(0,\tilde{r})& L_s^2&=C\cap \partial(B_-\cup B(z_k,s\varepsilon))\\
S_{\pm}&=\partial B_{\pm}& L^1&=L^2_\mu&
L^3&=C\cap\partial(B(0,\tilde{r})\cup B(z_k,\lambda\varepsilon))
\end{align*}
Here $\delta>0$ is a small number that will be fixed later, $s\in[\mu,\lambda]$ is a parameter.

\smallskip
\underline{\emph{Outline of the proof.}} The proof of the proposition will consist of the following steps.
\begin{itemize}
\item[\underline{Step 1.}] $L^1$ is a disconnected sum of $L_{-}$ and the link of singularity $\LS$;
\item[\underline{Step 2.}] $L^2_\lambda$ arises from $L^2_{\mu}$ by adding $p$ \ohands;
\item[\underline{Step 3.}] $L^3$ is isotopic to $L^2_\lambda$;
\item[\underline{Step 4.}] $L_+$ is isotopic to $L^3$.
\end{itemize}
The most important part is Step 2, all others are technical. The notation $L^1$, $L^2$ and $L^3$ suggests in which step does the given link appear.

\smallskip
\begin{figure}
\schempic
\begin{psfigure}\label{fig:schempic}
Schematic presentation of the proof of Proposition~\ref{singhand}. The curve $C$ (not drawn on the picture)
is intersected with boundaries of shaded sets providing
links $L_-$, $L^1$, $L^2$, $L^3$, and, finally, $L_+$.
\end{psfigure}
\end{figure}

In proving  Steps 2, 3 and 4 we will use the following lemma, which is a slight generalization of a standard result about isotopies. For
a convenience of the reader we present also a sketch of proof.

\begin{lemma}[Transverse isotopy]\label{tris}
Let $S^3=W_N\cup W_S$ be a decomposition of $S^3$ into an upper ''northern'' and lower ''southern'' closed hemispheres 
and let $S^2_{eq}=W_N\cup W_S$ be the ''equator''. We denote by $W_N^o$ and $W_S^o$ the interiors of $W_N$, respectively $W_S$. 
Assume that $\phi_s\colon S^3\to\mathbb{C}^2$ is a family of embeddings with
following assumptions.
\begin{itemize}
\item[(Is1)] $\phi\colon S^3\times[0,1]\to \mathbb{C}^2\times[0,1]$ given by $\phi(x,s)=(\phi_s(x),s)$ is continuous, i.e. $\phi_s$ is
a continuous family;
\item[(Is2)] $\phi_s$ is a smooth family when restricted to $W_N$ and to $W_S$, in particular it is smooth when restricted to $S^2_{eq}$;
\item[(Is3)] the image $\phi_s(W_N^o)$ and $\phi_s(W_S^o)$ is transverse to $C$;
\item[(Is4)] (the crucial in our applications) the image $\phi_s(S^2_{eq})$ is transverse to $C$.
\end{itemize}
Then the links $\phi_0^{-1}(C)$ and $\phi_1^{-1}(C)$ are isotopic.
\end{lemma}
\begin{proof}[Proof of Lemma~\ref{tris}]
If $\phi_s$ is $C^1$ smooth, the statement is standard. The proof in this case is slightly more technical, but follows the same pattern. Namely, we shall
prove that for any $s\in[0,1]$ and for any $s'$ sufficiently close to $s$, the links $\phi_s^{-1}(C)$ and $\phi_{s'}^{-1}(C)$ are isotopic and the statement 
shall follow from compactness and connectedness of the interval $[0,1]$.

Let us then consider a particular $s\in[0,1]$. Recall that $C$ was given by an equation $\{F=0\}$. Let $S^3_{reg}$, respectively $S^2_{eq,reg}$,
be the set of points
$x\in S^3$ (resp. $x\in S^2_{eq}$) such that $\phi_s(S^3)$ (resp. $\phi_s(S^2_{eq})$) is transverse to $F^{-1}(F(\phi_s(x)))$ at $\phi_s(x)$.

Now for each $x\in W_N\cap S^3_{reg}$ we can choose a vector $v_s^N(x)$ such that 
\begin{equation}\label{eq:vsdef}
DF\cdot\left(\frac{\partial\phi_s}{\partial s}+v_s^N(x)\right)=0
\end{equation}
(here $DF$ means the derivative regarded as a $4\times 2$ real matrix). 
This property means that $F\circ\phi_s$ is constant along the integral curves of the (non-autonomous) vector field $v_s^N$.
Now two different vectors $v_s^N(x)$ and $\tilde{v}_s^N(x)$ satisfying \eqref{eq:vsdef} differ by a vector which is tangent to $(F\circ \phi_s)^{-1}(F(\phi_s(x)))$.
In particular, we can pick $v_s^N(x)$ to be a smooth vector field, and, whenever $x\in S^2_{eq,reg}$, 
we can make $v_s^N(x)$ is tangent to $S^2_{eq}$. As 
each fiber $F^{-1}(F(\phi_s(x)))$ which is transverse to $S^2_{eq}$ intersects $S^2_{eq}$ in
finitely many points, we see that the vector fields $v_s^N$ is then uniquely defined on $S^2_{eq,reg}$.

Similarly we construct a vector field $v_s^S(x)$. The two vector vields $v_s^S$ and $v_s^N$ agree on $S^2_{eq,reg}$ and therefore they
can be glued to produce a vector field $v_s$ defined on $U=(S^3_{reg}\setminus S^2_{eq})\cup S^2_{eq,reg}$. As $v_s^S$ and $v_s^N$ are smooth, $v_s$
is locally Lipschitz. By Cauchy's theorem, $v_s$ can be integrated to a local diffeomorphism. This diffeomorphism 
maps fibers of $F\circ \phi_s$ to fibers of $F\circ \phi_{s'}$, for $s'$ sufficiently close to $s$.

Now the assumptions (Is3) and (Is4) guarantee that $\phi_s^{-1}(C)$ lies in the interior of $U$. Therefore, $\phi_s^{-1}(C)$ is isotopic
to $\phi_{s'}^{-1}(C)$ for $s'$ close to $s$ and we conclude the proof.
\end{proof}

Before we pass to the core of the proof of Proposition~\ref{singhand}, let us make an obvious, but important, remark. 
The order of tangency of each branch of $G_j$ of $C$ to $Z_j$ (see \eqref{eq:tangentlines})
is, by \eqref{eq:param}, $p_j\ge 2$. Therefore, a point $z\in C$ sufficiently close to $z_k$, the tangent space $T_zC$ is very close
to $Z_j$ for some $j$. In particular, if we can show transversality of some space $X\subset\mathbb{C}^2$ to all of $Z_j$, 
we can often claim the transversality of $X$ to $C$.

\smallskip
\noindent\underline{\emph{Step 1.}} By condition (S2) above, the intersection of $B_-$ and $B(z_k,\mu\varepsilon)$ is disjoint from $C$.
Therefore $C\cap (S_-\setminus B(z_k,\mu\varepsilon))=C\cap S_-=L_-$ and
$C\cap (S(z_k,\mu\varepsilon)\setminus B_-)=C\cap S(z_k,\mu\varepsilon)=\LS_k$. Thus the intersection of $C$ with
$\partial(B_-\cup B(z_k,\mu\varepsilon))$ is indeed a disjoint sum of $L_-$ and $\LS_k$.

\smallskip
\noindent\underline{\emph{Step 2.}} For any $s\in[\mu,\lambda]$, $C$ is transverse to $B(z_k,s\varepsilon)$ (because all $Z_j$'s are transverse
and $\varepsilon$ is sufficiently small). We are in a situation covered by Lemma~\ref{tris}: $\partial (B_-\cap B(z_k,s\varepsilon))$ can be regarded
as an image of a piecewise smooth map from $S^3$ to $\mathbb{C}^2$, which maps $S^3_S$ to $S_-$, $S^3_N$ to $S(z_k,\varepsilon)$
and $S^2_{eq}$ to $S_-\cap S(z_k,\varepsilon)$. Nevertheless, as the links $L^2_{\mu}$ and $L^2_{\lambda}$ 
are non-isotopic, some of the assumptions of Lemma~\ref{tris} must fail. Indeed, we shall show below that (Is4) is not satisfied 
(see Remark~\ref{ob:tangency} below) and we accomplish
Step~2 by studying the intersection of $C$ with $S_-\cap S(z_k,s\varepsilon)$.

\begin{figure}
\toymodel
\begin{psfigure}\label{fig:toymodel}
Toy model in three dimensions, which should help to understand Step~2.
Two balls $B_1$ and $B_2$. A plane $C$ intersects the boundary of $\partial(B_1\cup B_2)$ in two disjoint circles (left picture). 
If we push the ball $B_2$ inside $B_1$, this intersection becomes one circle. This is precisely a one handle attachment that occurs in Step 2.
\end{psfigure}
\end{figure}

Consider a branch $G_j$ of $C$ (see \eqref{eq:param}). The idea is that
up to terms of order $\tau^{p_j+1}$ or higher, the image of the branch $G_j$ is a $p_j$-times
covered disk, which lies in $Z_j$, so the situation described on Figure~\ref{fig:toymodel} happens precisely $p_j$ times, which gives $p_j$ \ohands.
Since the multiplicity of a singular point is equal to the sum of multiplicities of branches, this will conclude the proof.

\begin{figure}
\notationonAj
\begin{psfigure}
Schematic presentation of notation used in Step 2. The branch in question as multiplicity $p_j=3$. For cleareness of the picture, we draw only one disk $D_{ja}$
and do not label all objects. We also draw only a part of $\partial R_{s2}$, the whole $\partial R_{s2}$ is the full circle.
\end{psfigure}
\end{figure}

To be more rigorous, consider a disk
\[G_j \cap B(z_k,\lambda\varepsilon),\] 
which can be prezented as $\psi_j(R_\lambda)$, where $\psi_j$ is the parametrization of $G_j$ (see \eqref{eq:param}) and
\[
R_\lambda=\{\tau\in\mathbb{C}\colon (|\beta_j|^2+|\alpha_j|^2)|\tau|^{2 p_j}+\dots\le \lambda^2\varepsilon^2\},
\]
where $\dots$ denote higher order terms in $\tau$. Let
\[\Gamma=\psi_j^{-1}(B_-)\cap R_\lambda\] 
and for $s\in[\mu,\lambda]$, let
\[R_s=\psi_j^{-1}(B(z_k,\varepsilon s))\cap R_\lambda.\]
Observe that 
\begin{equation}\label{eq:psils}
\psi_j^{-1}(L^2_s)=\partial(\Gamma\cup R_s).
\end{equation}
It is also useful have in mind the following fact.

\begin{remark}\label{ob:tangency}
The intersection of the branch $G_j$ with $S_-\cap S(z_k,s\varepsilon)$ is not transverse (and so the condition (Is4) of Lemma~\ref{tris} is not satisfied, and
so one may expect a change of topology of link $L^{2}_s$) if and only if $\partial\Gamma$ is transverse to $\partial R_s$.
\end{remark}

Using the local parametrization, we can see that $R_s$, up to higher order terms, is given by 
\[|\tau|^2\le\left(\frac{s^2\varepsilon^2}{|\alpha_j|^2+|\beta_j|^2}\right)^{1/p_j}+\dots,\]
i.e. this is, up to higher order terms, a disk. In particular it is a convex set (see Remark~\ref{rem:c2close} below).
On the other hand we can compute explicitely the parametrization of $\partial\Gamma$. By plugging \eqref{eq:param} into 
the condition $|w_1|^2+|w_2|^2=(\rho_k-\varepsilon)^2$, and neglecting terms of order $p_j+1$ or higher in $\tau$ (and with $\varepsilon^2$), we get

\[
\partial\Gamma=\{\tau\colon \re \beta_j\tau^{p_j}=\frac12\varepsilon\rho_k\}.
\]

Chosing $\eta_j$ such that $\eta_j^{p_j}=\beta_j$, and writing in polar coordinates $(r,\phi)$ on~$R_\lambda$ 
\[\eta_j^{-1}\tau=r(\cos\phi+i\sin\phi)\] 
we finally obtain obtain

\begin{equation}\label{eq:gamma}
\partial\Gamma=\{(r,\phi)\in R_\lambda \colon r^{p_j}\cos p_j\phi=\frac12\varepsilon\rho_k\},
\end{equation}
modulo higher order terms.
We can see that $\partial\Gamma$ consists of $p_j$ connected components, indeed, for $\cos p_j\phi<0$
equation \eqref{eq:gamma} cannot hold. It follows that $\Gamma$ has also $p_j$ connected components, let us call them $\Gamma_{j1}\dots,\Gamma_{jp_j}$.
Each set $\Gamma_{ja}$ is convex. This follows from \eqref{eq:gamma} and a simple analytic observation, which we now state explicitely.

\begin{remark}\label{rem:c2close}
In general, the convexity of the connected 
subset of a disk given by $\{f\ge 0\}$ for some $f$ depends only on second derivatives of $f$. So if a function $g$ is $C^2$-close
enough to $f$, and the set $\{f\ge 0\}$ is convex, then $\{g\ge 0\}$ is convex, as well. Since the terms we neglect in the discussion above 
are of order $\tau^{p_i+1}$ and $p_i\ge 2$,
the convexity of $R_s$ follows from the convexity of a disk of radius $\varepsilon(|\alpha_j|^2+|\beta_j|^2)s^{1/p_j}$ and the convexity of $\Gamma_{ja}$
follows from the convexity of the set with bounaries parametrized by \eqref{eq:gamma} without higher order terms. Here we use implicitely condition (S4).
\end{remark}

Now consider a single $a=1,\dots,p_j$. By conditions (S2) and (S3) above $\Gamma_{ja}\cap R_{\mu}=\emptyset$ and
$\Gamma_{ja}\cap R_{\lambda}\neq\emptyset$.
Thus, by convexity, there exist a single $s=s_{ja}$ such that $\partial\Gamma_{ja}$ is tangent to $R_{s_{ja}}$. In particular, there are $p_j$ points on $R_\lambda$
such that $\partial\Gamma$ is tangent to $R_s$ for some $s$. Let us call them $y_{j1},\dots,y_{jp_j}$. Let us pick a very small disk $D_{ja}$ near $y_{ja}$.
Then for $s<s_{ja}$ close to $s_{ja}$, $\psi_j^{-1}(L^2_s)\cap D_{ja}$ (cf. \eqref{eq:psils}) consists of two arc: one on $\partial\Gamma$ and
the other on $\partial R_s$, see Figure~\ref{fig:nextfigure}. 
On the other hand, for $s>s_{ja}$ close to $s_{ja}$, $\psi_j^{-1}(L^2_s)\cap D_{ja}$ consists of two arc, each of them lies
partially on $\partial\Gamma$ and partially on $\partial R_s$. It follows that a \ohand{} addition occurs in $D_{ja}$ when $s$ passes through $s_{ja}$.

\begin{figure}
\nextpictureonC
\begin{psfigure}\label{fig:nextfigure}
Passing through $s_{ja}$. The picture presents $\psi_j^{-1}(B_-\cup B(z_k,s\varepsilon))\cap R_\lambda=\Gamma\cup R_s$, lying inside the disk $D_{ja}$.
On the left $s<s_{ja}$ and $\Gamma$ is disjoint from $R_s$, on the right $s>s_{ja}$ and $\Gamma\cap R_s\neq\emptyset$. The boundary of $\Gamma\cup R_s$
is mapped onto link $L_s^2$: we see that the topology changes by the \ohand{} addition as $s$ crosses $s_{ja}$. 
\end{psfigure}
\end{figure}
\medskip
\noindent\underline{\emph{Step 3.}} We isotope the ball $S_-=S(0,\rho_k-\varepsilon)$ to $S_+=S(0,\tilde{r})$ and use Lemma~\ref{tris}.
More precisely, consider a family of sets
\[B_s^3:=B(z_k,\lambda\varepsilon)\cup B(0,s),\]
where $s\in[\rho_k-\varepsilon,\tilde{r}]$. We can easily find a piecewise smooth family of maps $\phi_s^3\colon S^3\to\partial B_s^3$, such
that $\phi^3_s(W_N)\to S(0,s)$, $\phi^3_s(W_S)\to S(z_k,\lambda\varepsilon)$ and $\psi^3_s(S^2_{eq})=S(z_k\lambda\varepsilon)\cap S(0,s)$ (notation
from Lemma~\ref{tris}). Now $\phi^3_s(W_N^o)$ is transverse to $C$. Indeed, this follows by (S5) and the fact that  $z_k$ is not
in the image $\phi^3_s(W_N^o)$. Obviously $\phi^3_s(W_S^o)$ is tranverse to $C$, because $C$ is transverse to $S(z_k,\lambda\varepsilon)$. Therefore,
the condition (Is3) of Lemma~\ref{tris} is satisfied. We need to show (Is4). But observe that 
\begin{equation}\label{eq:whytilder}
S(0,\tilde{r})\cap S(z_k,\lambda\varepsilon)=S(z_k,\varepsilon)\cap\{w_2=0\}.
\end{equation}
Each tangent line $Z_j$ (see \eqref{eq:tangentlines}) is in fact transverse to $S(z_k,\lambda\varepsilon)\cup S(0,s)$ 
for all $s\in[\rho_k-\varepsilon,\tilde{r}]$.
(This follows from elementary geometric argument which we leave as an exercise. Figure~\ref{fig:step3} explains the key point of the argument, namely
that $\lambda$ has been chosen large enough.) Then by chosing $\varepsilon$ small enough we can ensure that $C$ is transverse to 
$S(z_k,\lambda\varepsilon)\cup S(0,s)$ so (Is4) is satisfied and the step is accomplished.
\begin{figure}
\whykappa
\begin{psfigure}\label{fig:step3}
Step 3. We explain, why the condition (S3) is important. $S_\mu$ is shorthand
for $S(z_k,\mu\varepsilon)$. The dotted ellipse represents $S_-\cap S_\lambda$. On the right hand side, there is one branch of $C$, namely $G_3$,
which doesn't intersect $S_-\cap S_\lambda$, if we start enlarging $S_-$, the intersection of $S_-\cap S_\lambda$ will eventually become
non-empty, so we shall meet a non-transversality point. If we choose $\lambda$ large enough, then all nontransversality points are dealt in with Step~2.
\end{psfigure}
\end{figure}

\medskip
\noindent\underline{\emph{Step 4.}} 
Let $B^4_0=B(0,\tilde{r})\cup B(z_k,\varepsilon)$. 
With a notation
of Lemma~\ref{tris}, let us consider a family of maps $\phi_s^4\colon S^3\to\mathbb{C}^2$ such that $\psi_s^4(W_N)=S_+\setminus B(z_k,\varepsilon)$
(in fact, we may assume that $\phi_s^4|_{W_N}$ does not depend on $s$), $\phi_0^4(W_S)=S(z_k,\varepsilon)\setminus B_+$ and $\phi_1^4(S^3)=S_+$.
Then the transversality of $\phi_s^4(W_N^o)$ and of $\phi_s^4(S^2_{eq})$ to $C$ (part of condition (Is3) and the condition (Is4)) is obvious).
It is not difficult to choose $\phi_s$ so that $\phi_s^4(W_S^o)$ is transverse to $C$. For example, one can observe that for any $s=[0,1]$,
the sphere
\[S_s=S\left(s\cdot z_k,\sqrt{(1-s)^2\rho_k^2+\lambda^2\varepsilon^2}\right)\]
passes through the intersection of $S(0,\tilde{r})\cap S(z_k,\lambda\varepsilon)$, for $s=0$ we have $S_0=S(z_k,\lambda\varepsilon)$ and for $s=1$, 
$S_1=S(0,\tilde{r})$. 
Then we can easily construct $\phi_s^4$ such that
$\phi_s^4(W_S)$ lies on $S_s$. It is a matter of direct computations to check that $\phi^4_s(W_S)$ is transverse to each tangent line $Z_j$
(see \eqref{eq:tangentlines}) so, if $\varepsilon$ is small enough, also to $C$. See Figure~\ref{fig:step4figure}.
\end{proof}
\begin{figure}
\stepfour
\begin{psfigure}\label{fig:step4figure}
Step 4. A schematic presentation of an isotopy of $\phi^4_s$. The consecutive images $\phi_s^4(W_S)$ are drawn with dashed lines, only
$\phi_0^4(W_S)$ and $\phi_1^4(W_S)$ (not labelled on the picture) are bold solid lines. 
The lines $Z_1$ and $Z_2$ are examples of possible tangent lines to $C$, they are all
transverse to images $\phi_s^4(W_S)$ for $s\in[0,1]$.
\end{psfigure}
\end{figure}

\begin{figure}
\begin{pspicture}(-8,-0.8)(3,1.5)
\rput{90}(-5,0){\psscalebox{0.25}{\includegraphics{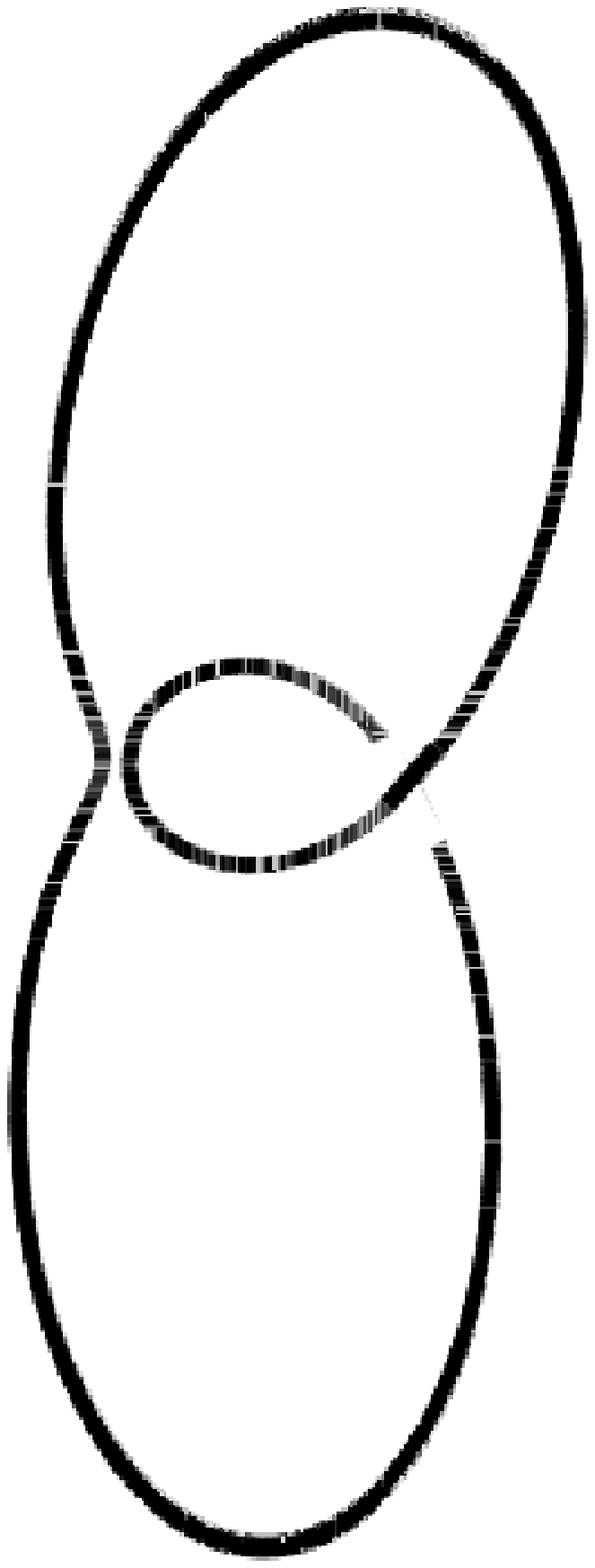}}}
\rput{90}(0.5,0){\psscalebox{0.25}{\includegraphics{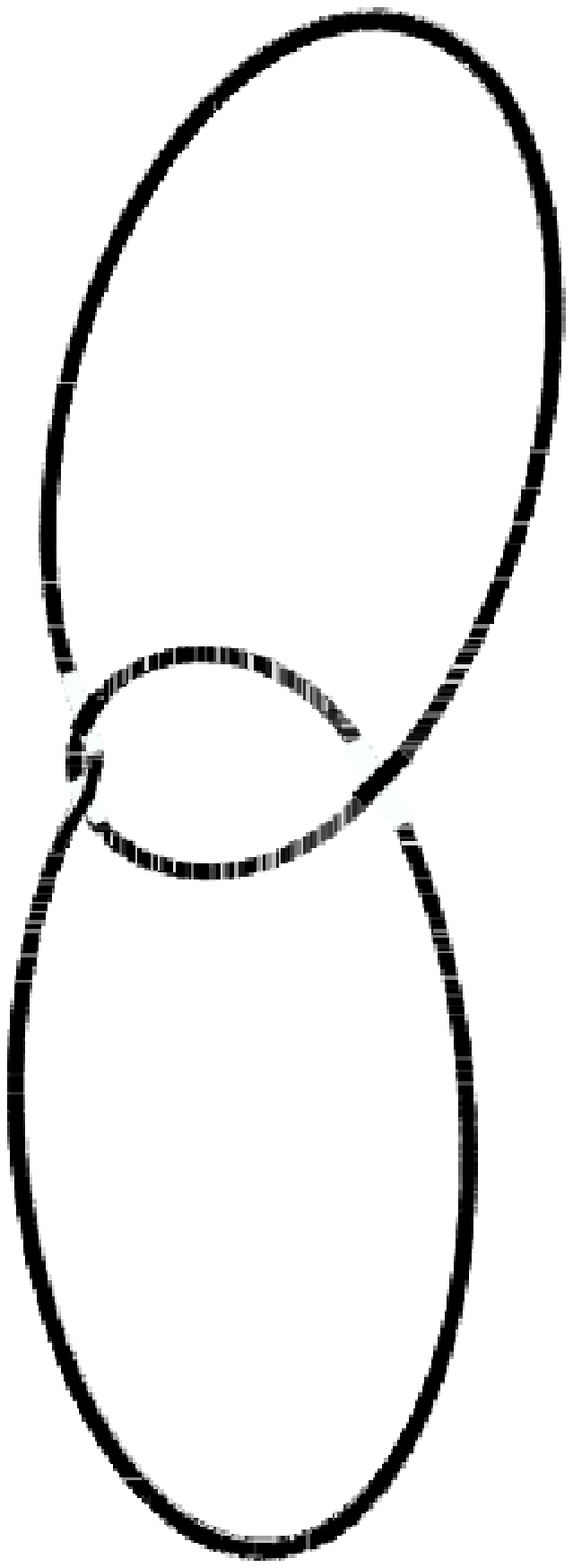}}}
\end{pspicture}
\begin{psfigure}\label{figA1}
Curve $\{x^3-x^2-y^2=0\}$ intersected with a sphere $S((-1,0),0.95)$ on the left and $S((-1,0),1.04)$ on the right. For radius $r=1$
we cross an ordinary double point. The trivial knot (on the left) becomes a trefoil after a change of one undercrossing to an 
overcrossing. 
(Figures \ref{figA1} and \ref{figA2} have been made using a C++ computer program written by the author. The author can provide the source code.)
\end{psfigure}
\end{figure}

Let us fix an arbitrary ordering of \ohands{}  at a given singular point once and for all. 
We shall then denote them $\tilde{H}_1,\dots,\tilde{H}_p$.
We can think of the procedure described 
in Proposition~\ref{singhand} as follows: first we take the disconnected sum of $L_-$ with $\LS$. After that we glue the handle $\tilde{H}_1$,
then $\tilde{H}_2$ and so on.
In this setting $\tilde{H}_1$ is a join handle 
and others
are either divorces or joins or marriages. Such handles will be called \emph{fake joins}, \emph{fake divorces} and \emph{fake
marriages} respectively. The total number of such handles at a point $z_k$ will be denoted $f^k_j$, $f^k_d$ and $f^k_m$. These numbers can
be computed by studying changes of the number of components and the Euler characterisics between  $C_-$ and $C_+$ and between $L_-$ and $L_+$ 
(see the proof of Proposition~\ref{prop:geom} below) and as such, they are independent of the ordering of handles.
\begin{example}
If $z_k$ is an ordinary double point (locally defined by $\{xy=0\}$), then $L_+$ arises from $L_-$ by changing
a negative crossing on some link diagram to a positive crossing (see Figure~\ref{figA1} and its explanation on the Figure~\ref{figA3}).
\end{example}

\begin{figure}
\begin{pspicture}(-8,-0.8)(3,2)
\rput(-5,0){\psscalebox{0.25}{\includegraphics{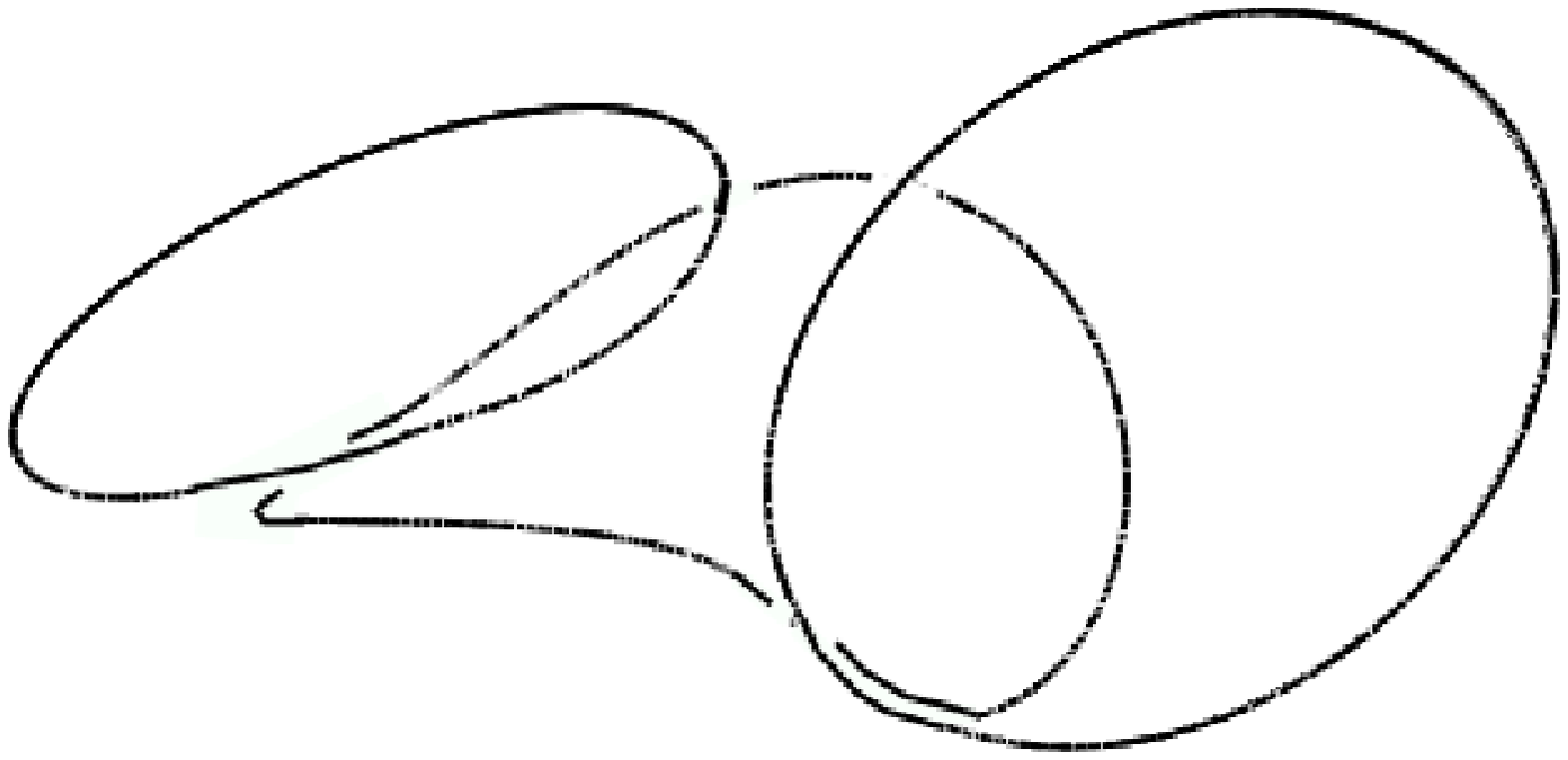}}}
\rput(0,0){\psscalebox{0.25}{\includegraphics{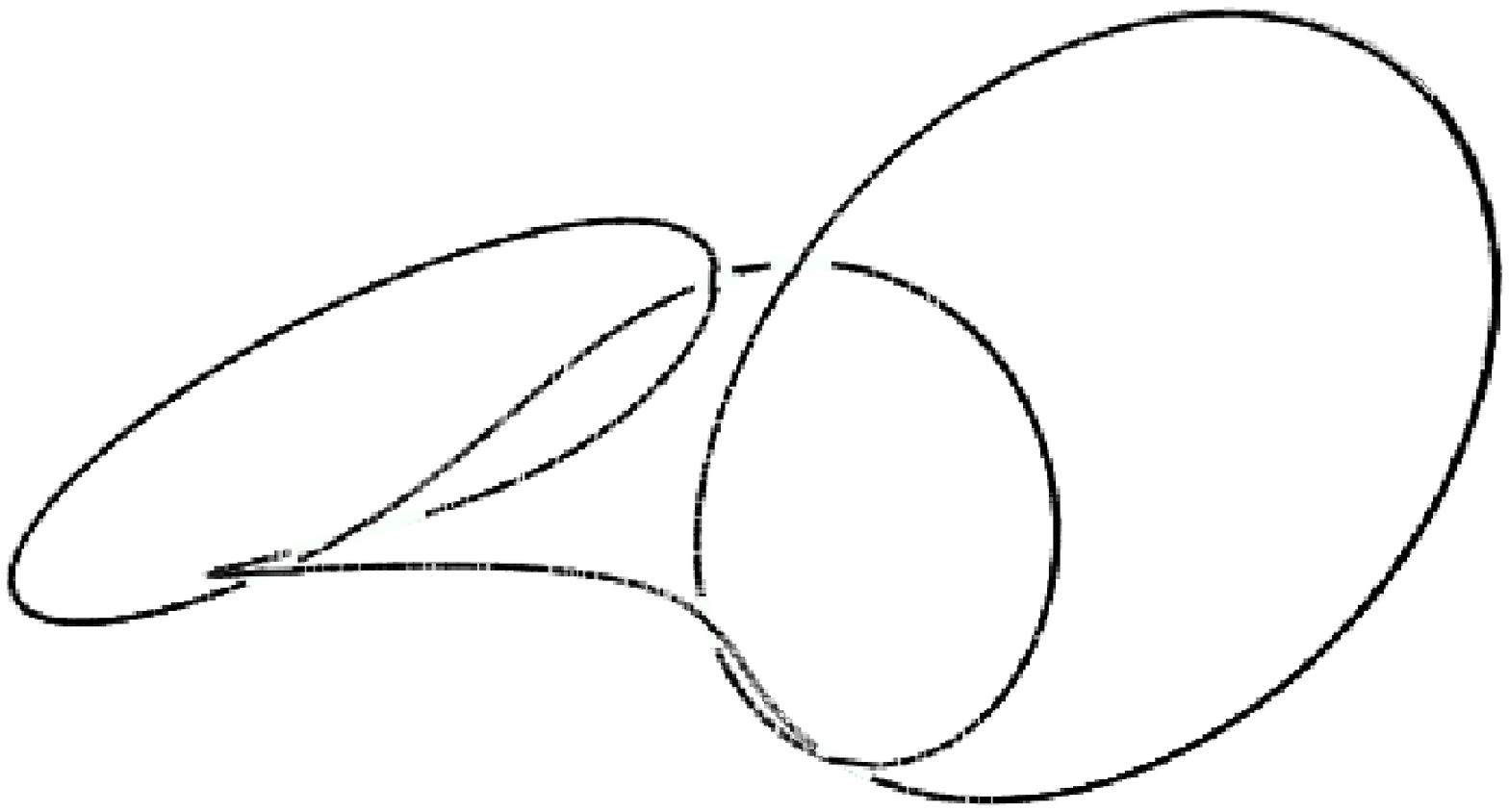}}}
\end{pspicture}
\begin{psfigure}\label{figA2}
Swallowtail curve (given in parametric form by $x(t)=t^3-3t$, $y(t)=t^4-2t^2$) 
intersected with a sphere $S((0,0),2.15)$ on the left and $S((0,0),2.5)$ on the right. 
We cross two $A_2$ singularities at $r=\sqrt{5}$. The two external circles on left picture twist around the middle one,
after crossing a singular point.
\end{psfigure}
\end{figure}

\begin{figure}
\gluepic
\begin{psfigure}\label{figA3}
The transformation of links shown on Figures~\ref{figA1} and \ref{figA2} explained as taking a sum with a Hopf link (resp. torus knot
$T_{2,3}$) and gluing two \ohands{} to the result. The bold parts of links represent places, where the handles are attached. Remark that
on Figure~\ref{figA2} the procedure is applied twice, because we cross two singular points at one time (i.e. $(0,0)$ violates the genericity
condition G1 in this case).
\end{psfigure}
\end{figure}
\section{Number of non--transversality points}\label{notran}
This section is auxillary in the sense that it provides some control over the number of non-transversality points, which might be useful in
the future. We only use one result from this section, namely the finiteness of critical points of an algebraic curve.

Let us consider a curve $C=\{F=0\}$ in $\CC^2$, such that $F$ is a reduced polynomial of 
degree $d$.
Let $\xi=(\xi_1,\xi_2)\in\CC^2$ be a fixed point (a ball centre). Let $S_r=S(\xi,r)$ be a 
three--sphere of radius $r$
centered at $\xi$. Let $w=(w_1,w_2)$ be an arbitrary point in $C\cap S_r$. Assume that $C$ is smooth at $w$.

\begin{lemma}\label{lem3.1}
The intersection $C\cap S_r$ is transverse at $w$ if and only if the determinant
\[
J_\xi(w)=\det\left(\begin{matrix}
\overline{\frac{\partial F}{\partial w_1}}(w)&\overline{\frac{\partial F}{\partial w_2}}(w)\\
w_1-\xi_1&w_2-\xi_2
\end{matrix}\right)
\]
does not vanish.
\end{lemma}
\begin{proof}
Assume that $C$ is not transverse to $S_r$ at $w$. This means that
\[T_wC+T_wS_r\neq \CC^2.\]
Since $T_wS_r$ is real three dimensional, $T_wC+T_wS_r=T_wS_r$, thus
\[T_wC\subset T_wS_r.\]
Taking the orthogonal complements of these spaces we see that
\[N_wS_r\subset N_wC.\]
But $N_wC$ is a complex space. Thus $\mathbf{i}\cdot N_wS_r\subset N_wC$ and by dimension arguments 
we get that
\[N_wS_r\otimes\CC=N_wC.\]
Now $N_wS_r\otimes\CC$ is spanned over $\CC$ by a vector $(w_1-\xi_1,w_2-\xi_2)$.
The lemma follows (the above reasoning can be reversed to show the ''if'' part). 
\end{proof}
If $w$ is a singular point of $C$, $J_\xi(w)=0$ by the definition.
\begin{corollary}
For a curve $C$ of degree $d$ and a generic point $\xi\in\CC^2$ there are $d(d-2)$ such points 
$($counted with multiplicities$)$
$w\in C$ where the intersection
\[C\cap S(\xi, ||w-\xi||),\]
is not transverse at $w$.
\end{corollary}
\begin{proof}
For a fixed $\xi$, $J_\xi(w)$ is a polynomial of degree $d-1$ in $w$ and $1$ in $\bar{w}$. 
Intersecting $\{J_\xi=0\}$ with $C$ of degree $d$ yields $d^2-2d$ points (counted with multiplicities)
by generalised B\'ezout theorem (see e.g. \cite[Theorem 1]{Chen}). 
\end{proof}
\begin{remark}\label{rem:fin}
The number of intersection points can be effectively larger than $d^2-2d$: as the 
curve $\{J_\xi=0\}$
is not complex, there might occur intersection points of multiplicity $-1$. Anyway, this number is always finite, because both $C$ and $J_\xi$ are
real algebraic.
\end{remark}

The local intersection index of $C$ with $\{J_\xi(w)=0\}$ at a singular 
point $z$ can be effectively calculated. We have the following lemma.
\begin{lemma}
Assume that $0\in\mathbb{C}^2$ is a singular point of $C$.
The local intersection index of $C$ with $\{J_\xi=0\}$ at $0$
is equal to the Milnor number $\mu$ of $C$
at $0$ minus $1$.
\end{lemma}
\begin{proof}
This follows from Teissier lemma (see \cite{Pl} or \cite{GP}), which states that
\[(f,J(f,g))_0=\mu(f)+(f,g)_0-1,\]
where $(a,b)_0$ denotes the local intersection index of curves $\{a=0\}$ and $\{b=0\}$ at $0$
and $J(f,g)$ is the Jacobian
\[\frac{\partial f}{\partial w_1}\frac{\partial g}{\partial w_2}-\frac{\partial f}{\partial w_2}\frac{\partial g}{\partial w_1}.\]
We shall apply this lemma to the case when $f=F$ is the polynomial defining the curve $C$, whereas
$g$ is the distance function:
\[g(w_1,w_2)=|w_1-\xi_1|^2+|w_2-\xi_2|^2\]
Then $(f,g)_0=0$. 
In fact, intersection of $\{f=0\}$ and $\{g=0\}$ is real one dimensional. 
But if we perturb $g$ to $g-\mathbf{i}\varepsilon$
the intersection set becomes empty.

The issue is that the Teissier lemma holds when $f$ and $g$ are holomorphic. To see that nothing
bad happens, if $g$ is as above, we have to skim 
through a part of the proof of Teissier lemma (see e.g. \cite{Pl}).
Assume for a while that the curve $\{f=0\}$ can be parametrised
near $0$ by
\[w_1=t^n,\,\,\,w_2=w_2(t),\]
where $w_2(t)$ is holomorphic and $n$ is the local multiplicity of $\{f=0\}$ at 0. 
(The case of many branches does not present new difficulties.)
Then
\begin{equation}\label{eq:teissier}
\begin{split}
\frac{\partial f}{\partial w_1}(t^n,w_2(t))\cdot nt^{n-1}+\frac{\partial f}{\partial w_2}(t^n,w_2(t))&=0\\
\frac{\partial g}{\partial w_1}(t^n,w_2(t))\cdot nt^{n-1}+\frac{\partial g}{\partial w_2}(t^n,w_2(t))&=\frac{d}{dt}g(t^n,w_2(t)).
\end{split}
\end{equation}
The first equation follows from differentiating the identity $f(t^n,w_2(t))\equiv 0$. The second is simply the chain
rule applied to its r.h.s. On its l.h.s. we 
could have terms with
$\frac{\partial g}{\partial\bar{w}_2}\frac{\partial\bar{w_2}}{\partial t}$. But they vanish, as $w_2$ is holomorphic.

From \eqref{eq:teissier} we get
\begin{equation}\label{eq:ordercomp}
nt^{n-1}J(f,g)(t^n,w_2(t))=-\frac{dg(t^n,w_2(t))}{dt}\cdot\frac{\partial f}{\partial w_2}(t^n,w_2(t)).
\end{equation}
Now we can compare orders with respect to $t$.
On the l.h.s. of \eqref{eq:ordercomp} we have
\[(n-1)+(f,J(f,g))_0.\]
Whereas on the r.h.s. we get
\[(f,g)_0-1+(f,\frac{\partial f}{\partial w_2})_0,\]
And we use another lemma, due also to Teissier, that $(f,\frac{\partial f}{\partial w_2})_0=\mu(f)+n-1$.
This can be done directly as $f$ is holomorphic.
\end{proof}

\section{Signature of a link and its properties}
Let $L\subset S^3$ be a link and $V$ a Seifert matrix of $L$ (see e.g. \cite{Ka} for necessary definitions). 
\begin{definition}\label{defsig}
Let us consider the symmetric form
\begin{equation}\label{eq:SF}
V+V^T.
\end{equation}
The \emph{signature} $\sigma(L)$ of $L$ is the signature of the above form. The
\emph{nullity} (denoted $n(L)$) is $1$ plus the dimension of a maximal null-space of the form \eqref{eq:SF}.
\end{definition}

The signature is an important knot cobordism invariant. 
Unlike many other invariants, signature behaves well under a \ohand{} addition.
More precisely we have
\begin{lemma}\label{mur} \emph{(}see \cite{Mur}\emph{)}
\begin{itemize}
\item[(a)] Let $L$ and $L'$ be two links such that $L'$ can be obtained from $L$ by a hyperbolic transformation
(see Definition~\ref{hyper} above). Then
\begin{align*}
 |n(L)-n(L')|&=1 \text{ and }\sigma(L)=\sigma(L')\text{; \ or}\\
 |\sigma(L)-\sigma(L')|&=1\text{ and }n(L)=n(L').
\end{align*}
\item[(b)] Signature is additive under the connected sum. 
The nullity of a connected sum of links
$L_1$ and $L_2$ is equal to $n(L_1)+n(L_2)-1$.
\item[(c)] Let $L$ be a link and $L'$ be a link resulting in the change from an 
undercrossing to an overcrossing on some planar diagram of $L$. Then either
\begin{align*}
&\sigma(L')-\sigma(L)\in\{0,-2\}\text{ and }n(L)=n(L')\text{; \ or}\\
&\sigma(L')=\sigma(L)-1\text{ and }|n(L)-n(L')|=1.
\end{align*}
\item[(d)]  $n$ does not exceed the number of components of the link.
\item[(e)] The signature and nullity are additive under the disconnected sum. 
\end{itemize}
\end{lemma}

The signature of a torus knot was computed for example in \cite{Ka,Li}.

\begin{lemma}\label{torsig}
Let $p,q>1$ be coprime numbers and $T_{p,q}$ be the $(p,q)$-torus knot.
Let us consider a set
\[\Sigma=\left\{\frac{i}{p}+\frac{j}{q},1\le i<p, 1\le j<q\right\},\]
\emph{(}note in passing that this is the spectrum of the singularity $x^p-y^q=0$, see \cite{BN} for a detailed discussion of this phenomenon\emph{).} Then 
\begin{equation}
\sigma(T_{p,q})=\#\Sigma-2\#\Sigma\cap (1/2,3/2).
\end{equation}
\end{lemma}
This means that $\sigma$ counts the elements in $\Sigma$ with a sign $-1$ or $+1$ according
to whether the element lies in $(1/2,3/2)$ or not.
\begin{example} We have
\begin{equation}\label{eq:smallsig}
\begin{split}
\sigma(T_{2,2n+1})&=-2n;\\
\sigma(T_{3,n})&=4\intfrac{n}{6}-2(n-1);\\
\sigma(T_{4,n})&=4\intfrac{n}{4}-3(n-1).
\end{split}
\end{equation}
Moreover, for $p$ and $q$ large, $\sigma(T_{p,q})=-\frac{pq}{2}+\dots$, where $\dots$ denote lower order terms in $p$ and $q$.
\end{example}
Lemma~\ref{torsig} holds even if $p$ and $q$ are not coprime (see \cite{Ka}): then we have a torus link instead of a knot.

Next result is a direct consequence of the discussion in \cite{Neu0}. It holds, in fact, for any graph link with non-vanishing
Alexander polynomial.
\begin{lemma}\label{lem:folklore}
Let $L$ be an algebraic link. Then $n(L)=c(L)$.
\end{lemma}

The following result of A. N\'emethi \cite{Nem2} will also be useful
\begin{proposition}\label{nemeth1}
Let $f$ be a reduced polynomial in two variables such that the curve $\{f=0\}$ has an isolated
singularity at $(0,0)$. Let $f=f_1\cdot f_2$ be the decomposition of $f$ locally near $(0,0)$, such
that $f_1(0,0)=f_2(0,0)=0$. Let $L$, $L_1$ and $L_2$ be the links of singularities of $\{f=0\}$, $\{f_1=0\}$
and $\{f_2=0\}$ at $(0,0)$ and $\sigma$, $\sigma_1$, $\sigma_2$ its signatures.
Then we have
\[\sigma\le \sigma_1+\sigma_2.\]
\end{proposition}
We could use the proof from~\cite{Nem1}. Nevertheless, we shall show a topological proof at the end of next section.

\begin{lemma}\label{sigmult}
Let $L$ be a link of plane curve singularity with $r$ branches. Then $\sigma(L)\le 1-r$. 
Moreover the equality holds only for the Hopf link and a trivial knot.
\end{lemma} 
\begin{proof}
Let $G$ be a germ of a singular curve bounding $L$. Let $\mu$ be the Milnor number of the singularity of $G$ and $\delta=\frac12(\mu+r-1)$
be the $\delta-$invariant of the singular point. There is a classical result (see e.g. \cite{Nem3}) that $-\sigma(L)\ge\delta$. This
settles the case if $r=1$. If $r>2$ we use the inequality $\delta\ge\frac12r(r-1)>r$ (which holds because $2\delta\ge \sum_{i\neq j}(C_i\cdot C_j)$, 
where $(C_i\cdot C_j)$ is the intersection index of two branches at a given singular point) and we are done. If $r=2$ we know that
$\delta\ge 1$, with equality only for an ordinary double point.
\end{proof}
\begin{corollary}\label{stupid}
Let $L=K_1\cup\dots\cup K_{n+1}$ be a link of a plane curve singularity with $n+1$ branches. Then
\[\sigma(L)\le \sigma(K_{n+1})+1-n.\]
\end{corollary}
\begin{proof}
Let $L'=K_1\cup\dots\cup K_n$. By Proposition~\ref{nemeth1} $\sigma(L)\le \sigma(L')+\sigma(K_{n+1})$. By Lemma~\ref{sigmult},
$\sigma(L')\le 1-n$.
\end{proof}

\section{Changes of signature upon an addition of a handle}\label{classsig}
In order to study the behaviour of some invariants of knots
let us introduce the following notation. Here $r\in\mathbb{R}$, $r>0$ and $r\not\in\{\rho_1,\dots,\rho_n\}$.
\begin{itemize}
\item $L_r$ the link $C\cap S(\xi,r)$;
\item $C_r$ the surface $C\cap B(\xi,r)$ and $\hat{C}_r$ is its normalization;
\item $k(C_r)$ number of connected components of $\hat{C}_r$;
\item $c(C_r)$ or $c(L_r)$ number of boundary components of $C_r$;
\item $\chi(C_r)$ the Euler characteristic of $C_r$;
\item $\pg(C_r)$ the genus of $C_r$, which for smooth $C_r$ satisfies $2k-2\pg=\chi+c$;
\item $\sigma(L_r)$ the signature of $L_r$
\item $n(L_r)$ the nullity of $L_r$.
\end{itemize}
If $C_r$ is singular, we are interested in the geometric genus
of $C_r$, i.e. the genus of normalisation of $C_r$. This explains the notation $\pg$ for a genus.

The following table describes the change of the above quantities upon attaching a handle.

\medskip
\begin{center}
\begin{tabular}{|p{2cm}|c|c|c|c|c|c|c|c|}\hline\label{tab}
name&index&$\Delta c$&$\Delta k$&$\Delta\chi$&$\Delta \pg$&$\Delta\sigma$&$\Delta n$\\\hline
birth&0&1&1&1&0&0&1\\\hline
death&2&-1&0&1&0&0&-1\\\hline
join&1&-1&-1&-1&0&$s$&$s'$\\\hline
divorce&1&1&0&-1&0&$s$&$s'$\\\hline
marriage&1&-1&0&-1&1&$s$&$s'$\\\hline
\end{tabular}
\end{center}
\medskip
Here $s,s'\in\{-1,0,1\}$ and $|s|+|s'|=1$ by Lemma~\ref{mur}~(a). 

\smallskip
Let
\begin{equation}\label{eq:w}
\begin{split}
w(L)&=-\sigma(L)+n(L)-c(L)\\
u(L)&=-\sigma(L)-n(L)+c(L)
\end{split}
\end{equation}

\begin{lemma}
If $L$ is a non-trivial link of singularity then $u(L)>0$ and $w(L)\ge 0$. Moreover, $w(L)=0$
if and only if $L$ is a Hopf link.
\end{lemma}
\begin{proof}
We use Lemma~\ref{sigmult} to prove this for $w(L)$. For $u(L)$ we use the fact that the signature is negative and Lemma~\ref{mur}(d).
\end{proof}

For a knot, by Lemma~\ref{mur}(d) we have $w(L)=u(L)=-\sigma(L)$. In general case of links we have
\begin{equation}\label{eq:boundwu}
\begin{split}
-\sigma(L)+(c(L)-1)\ge u(L)\ge -\sigma(L)\ge \\
\ge w(L)\ge -\sigma(L)-(c(L)-1).
\end{split}
\end{equation}
\begin{lemma}\label{addit}
The invariants $w(L)$ and $u(L)$ are additive under the disconnected sum.\qed
\end{lemma}
\begin{lemma}\label{nodecrease}
Attaching a birth, death, marriage or join handle does not decrease $w(L)$. 
\end{lemma}
\begin{proof}
Only the case of \ohands\ requires some attention. The number of component decreases by 1 and
either the nullity or the signature can change, and only by 1.
\end{proof}
\begin{remark}\label{decrease}
The divorce handle might decrease the quantity $w(L)$ at most by $2$.
\end{remark}
\begin{lemma}\label{noincreaseu}
Attaching a birth, death, marriage or join handle does not increase $u(L)$. The divorce might increase $u(L)$ at most by $2$. \qed 
\end{lemma}
\begin{lemma}\label{singbound}
Let $z_k$ be a singular point of $C$, $\LS_k$ the link of its singularity and $f_d^k$ the number
of fake divorces (see comment after the proof of Proposition~\ref{singhand}) at $z_k$. Let, for $\varepsilon>0$ small enough $L_\pm=L_{\rho_k\pm\varepsilon}$, where
$\rho_k=||z_k-\xi||$. Then
\begin{align*}
w(L_+)&\ge w(L_-)+w(\LS_k)-2f_d^k\\
u(L_+)&\le u(L_-)+u(\LS_k)+2f_d^k
\end{align*}
\end{lemma}
\begin{proof}
We use the notation from the proof of Proposition~\ref{singhand}. We have
\begin{align*}
w(L^1)&=w(L_-)+w(\LS_k)&&\text{step 1}\\
w(L^2)&\ge w(L^1)-2f_d^k&&\text{step 2}\\
w(L_+)&=w(L^2)&&\text{steps 3 and 4.}
\end{align*}
In the middle equations we have used the fact that a fake divorce can lower the invariant at most by 2.
The proof for $u$ is identical.
\end{proof}
\begin{lemma}\label{boundhand}
 Assume that $C$ is smooth. Let $\pg$ be the
genus of the curve $C$ and $d$ the number of its components at infinity. Let also
$a_b$, $a_m$, $a_d$, and $a_j$ denote the number of birth, marriage, divorce and  join handles.
The following formulae hold
\begin{equation}\label{eq:nohandle}
\begin{split}
a_m&=\pg\\
a_b+a_d-a_j-a_m&=d\\
a_b-a_j&=1.
\end{split}
\end{equation}
In particular
\begin{equation}\label{eq:ad}
a_d=d+\pg-1.
\end{equation}
\end{lemma}
\begin{proof}
For $r<\rho_1$, $L_r$ is empty. Thus the first handle must be a birth and for $r\in (\rho_1,\rho_2)$, $L_r$ is an unknot. It has
$\pg=0$, $c=1$ and $k=1$. When we cross next critical points, these quantities change according to the table on page~\pageref{tab}.
For $r>\rho_n$ we have the link at infinity and $C_r$ is isotopic to $C$.
\end{proof}

\begin{proposition}\label{prop:geom}
Let $C$ be an algebraic curve in $\mathbb{C}^2$, not necessarily smooth.
For a generic point $\xi$, let $S_0=S(\xi,r_0)$ and $S_1=S(\xi,r_1)$ (with $r_0<r_1$) be two spheres intersecting transversally with $C$.
For $i=0,1$ we define $\pg_i=\pg(C_{r_i})$, $c_i=c(C_{r_i})$ and $k_i=k(C_{r_i})$.

Let $a_d^{01}$ and $f_d^{01}$ be the numbers of divorces, respectively fake divorces, 
on $C$, which lie between $S_0$ and $S_1$.
Then
\[a_d^{01}+f_d^{01}\le \pg_1-\pg_0+c_1-c_0-(k_1-k_0).\]
\end{proposition}
\begin{proof}
Let $\pi:\hat{C}\to C$ be the normalisation map. 
The composition of $\pi$ with the distance function
$g$ (see \eqref{eq:dist})
restricted to $C$ yields a function $\hat{g}:\hat{C}\to\mathbb{R}$. 
This function does not have to be a Morse function on $\hat{C}$,
but we can take a small subharmonic perturbation of $\hat{g}$ on $\hat{C_{r_1}}$, such that the resulting
function is Morse in the preimage $\pi^{-1}B(\xi,r_1)$. This perturbation we shall still denote by $\hat{g}$.
Let $\hat{a}_b$, $\hat{a}_d$, $\hat{a}_j$ and $\hat{a}_m$ be
the number of births, divorces, joins and marriages of $\hat{g}$ in $U=\pi^{-1}(B(\xi,r_1)\setminus B(\xi,r_0))$. We need the following result
\begin{lemma} There is a bound
\begin{equation}\label{eq:hatad}
\hat{a}_d\ge a_d^{01}+f_d^{01}.
\end{equation}
\end{lemma}
\begin{proof}
If $z_k\in C$ is a smooth point of $C$ and critical point of $g$ then $\pi^{-1}(z_k)$ is 
a critical point of $\hat{g}$ of the same index. Moreover, if $z_k$ is a divorce, join or marriage then
$\pi^{-1}(z_k)$ will also be, respectively, a divorce, join or a marriage. 

Next we show that any fake divorce on $C$ corresponds to a divorce on $\hat{C}$. 
This is done by comparing the changes of topology when crossing a singular point 
with the changes of topology of normalisation. So let $z_k$ be a singular point of $C$.
Let us define
\[C_\pm=C\cap B(\xi,\rho_k\pm\varepsilon)\text{ and } L_\pm=\partial C_\pm\]
Let $\hat{C}_\pm$ be the normalization. Define also
\[
\Delta_g=\pg(C_+)-\pg(C_-),\ \ \ 
\Delta_k=k(C_+)-k(C_-),\ \ \ 
\Delta_c=c(L_+)-c(L_-).
\]
Observe that from a topological (as opposed to smooth) point of view, passing through a singular point of multiplicity $p$ and $r$ branches
amounts to picking $r$ disks and 
attaching them to $\hat{C}_-$ with $p$ \ohands{}.
Analogously to \eqref{eq:nohandle} we get then $f^k_m=\Delta_g$, $f^k_d-f^k_j-f^k_m=\Delta_c$ and $f^k_j=\Delta_k$. Hence
\[f^k_d=\Delta_c+\Delta_g-\Delta_k\]
The number of divorces on $\hat{C}$ that are close to $\pi^{-1}(z_k)$  (denote this number by $\hat{a}^k_d$) can be computed in the same way.
Since the number of boundary components of $\hat{C}_\pm$ is the same as $c(C_\pm)$, and $\Delta_g$ measures also the change of genus between
$\hat{C}_+$ and $\hat{C}_-$, we have
\[\hat{a}^k_d=\Delta_c+\Delta_g-\Delta_k=f^k_d.\]
\end{proof}

\noindent\emph{Finishing the proof of Proposition~\ref{prop:geom}.}
Let us consider the changes of the topology of $\hat{C}\cap\hat{g}^{-1}((-\infty,r^2))$ as $r$ changes from $r_0$ to $r_1$. 
The number of components of the boundary changes by $c_1-c_0$, while the genus by $g_1-g_0$ and the number of connected components of 
normalization by $k_1-k_0$.
Using the table on page~\pageref{tab} (compare the argument in the proof of Lemma~\ref{boundhand}) 
we get $\hat{a}_d=g_1-g_0+c_1-c_0-(k_1-k_0)$.
\end{proof}

\begin{remark}
In most applications we will have $k_0=k_1=1$, for example in the case when $L_1$ is a link at infinity of a reduced curve and $L_0$ is trivial knot.
\end{remark}
\begin{example}\label{ex:conic}
Let $C$ be a curve given by $x^3-x^2-y^2=0$ (see Figure~\ref{figA1} above, but now the center is in a different place), 
$\xi=(0,0)$, $r_0$ is small and let us take $r_1$ large enough. Then
$L_0$ is the Hopf link, $L_1$ is the treefoil, $\pg_1=\pg_0=0$ ($C$ is rational), $c_0=2$, $c_1=1$, $k_1=1$ but $k_0=2$ ($\hat{C}_0$ consists of two
disks). Then the number of divorces is bounded by $0$ and indeed, there is only one critical value between $r_0$ and $r_1$ and the corresponding handle
is a join.
\end{example}

\begin{corollary}\label{cor:geom}
If $C\subset\mathbb{C}^2$ is a reduced plane algebraic curve and its link at infinity has $d$ components, then for any generic $\xi$ the total number
of divorces on $C$ (including the fake divorces) satisfies
\[a_d+f_d\le \pg(C)+d-1.\]
\end{corollary}
\begin{proof}
Let us pick a generic $\xi$ and choose $r_0\in (\rho_1,\rho_2)$ while $r_1$ is sufficiently large. Then $S_0$ is an unknot, because the first handle that occurs
when coming from $r=0$, is always a birth. Moreover, $S_1\cap C$ is the link of $C$ at infinity and so it has $d$ components. The statement follows
from Proposition~\ref{prop:geom}
\end{proof}

\begin{theorem}\label{mthm2}
Let $C$ be a curve with link at infinity $L_\infty$ and with singular points $z_1,\dots,z_n$, such that the
link at the singular point $z_k$ is $\LS_k$. Then
\begin{align*}
w(L_\infty)&\ge\sum_{k=1}^n w(\LS_k)-2(p_g(C)+d-1),\\
u(L_\infty)&\le\sum_{k=1}^n u(\LS_k)+2(p_g(C)+d-1),
\end{align*}
where $d$ is the number of components of $L_\infty$.
\end{theorem}
\begin{proof}
The proof now is straightforward. Let us take a generic $\xi$. Then, for $r\in(\rho_1,\rho_2)$, $L_r$ is an unknot (see Remark~\ref{rem:smooth2}), 
so $w(L_r)=u(L_r)=0$. Then, as
we cross subsequent singular points, $w(L_r)$ and $u(L_r)$ change (see Lemmas \ref{nodecrease}, \ref{decrease}, \ref{noincreaseu}
and \ref{singbound}). We obtain
\[w(L_\infty)\ge\sum_{k=1}^n (w(\LS_k)-2f_d^k)-2a_d\]
and similar expression for $u$. The theorem follows now from Corollary~\ref{cor:geom}. 
\end{proof}
\begin{remark}\label{rem:blind}
Observe that the first inequality in Theorem~\ref{mthm2} (as applications below show, the more important one) 
'does not see' ordinary double points, because if $z_k$ is an ordinary
double point then $w(\LS_k)=0$ (however $u(\LS_k)=2$). 
\end{remark}

As the whole discussion leading to Theorem~\ref{mthm2} was quite involved, we present some examples.
\begin{example}
Consider a curve $\{x^3-x^2-y^2=0\}$, see Example~\ref{ex:conic}. An ordinary double point at $(0,0)$ is the only singular point (it has $w_L=0$
and $u_L=2$). The link at infinity is a trefoil with $w=u=2$. The geometric genus of a curve is equal to $0$.
\end{example}

\begin{example}
Let $C$ be a swallowtail curve as in Figure~\ref{figA2}. It has two ordinary cusps (the corresponding links of singularities are trefoils)
and one ordinary double point, its geometric genus is $0$
and the link at infinity is the torus knot $T_{3,4}$, with $w=u=6$. The inequalities in Theorem~\ref{mthm2} read $6\ge 4$ (the first one) and $6\le 6$ 
(the second one).
\end{example}

\begin{example}
Consider a curve parametrised by $x(t)=t^4$, $y(t)=t^6+t^9$. It has a singular point at $(0,0)$. According to \cite{EN} the link of
this singularity (let us call it $L_1$) is a $(15,2)$ cable on the trefoil. The curve has also three other ordinary double points 
(corresponding to $t=\sqrt[3]{1+i}$, 
which can be found by solving the equations
$x(t)=x(s)$, $y(t)=y(s)$, $t\neq s$). The link at infinity $L_{inf}$ (see \cite{Neu2}) is a $(9,2)$ cable on the trefoil. According to Lemma~\ref{cable}
$w(L_{inf})-w(L_1)=\sigma(T_{15,2})-\sigma(T_{9,2})=6$ and the same formula holds for $u$. Theorem~\ref{mthm2} holds, because $6\ge 6+3\cdot 0$ (inequality
for $w(L)$) and $6\le 6+3\cdot 2$. 
\end{example}

A good number of possible examples can be found also in \cite{BZ0,BZ1}, where a detailed list of plane algebraic curves
with the first Betti number $1$ is presented and for each curve on the list, its signularities are given explicitely. 
We provide one example (point (w) in the list of \cite{BZ1}), where a divorce handle occurs.
\begin{example}
Consider a curve parametrized by $x(t)=t^2-2t^{-1}$, $y(t)=2t-t^{-2}$. It has three ordinary cusps and no other singularities. It follows that 
$\sum w(\LS_k)=\sum u(\LS_k)=6$.
The curve has two branches at infinity, corresponding to $t\to\infty$ and $t\to 0$. Each branch is smooth at infinity and 
tangent to the line at infinity with the tangency order $2$.
An application of the algorithm of \cite{Neu2} shows that the link at infinity can be represented by the following splice diagram.

\begin{pspicture}(-5,-1.7)(5,0.7)
\psline{<->}(-3,0)(3,0)\psline(-1.5,0)(-1.5,-1.5)\psline(1.5,0)(1.5,-1.5)
\pscircle[fillstyle=solid,fillcolor=black](-1.5,0){0.05}
\pscircle[fillstyle=solid,fillcolor=black](1.5,0){0.05}
\pscircle[fillstyle=solid,fillcolor=black](-1.5,-1.5){0.05}
\pscircle[fillstyle=solid,fillcolor=black](1.5,-1.5){0.05}
\pscircle[fillstyle=solid,fillcolor=black](0,0){0.05}
\rput(0,0.2){\psscalebox{0.8}{root}}
\rput(1.2,0.2){\psscalebox{0.8}{$1$}}
\rput(1.8,0.2){\psscalebox{0.8}{$1$}}
\rput(1.7,-0.3){\psscalebox{0.8}{$2$}}
\rput(-1.2,0.2){\psscalebox{0.8}{$1$}}
\rput(-1.8,0.2){\psscalebox{0.8}{$1$}}
\rput(-1.3,-0.3){\psscalebox{0.8}{$2$}}
\end{pspicture}

Then, the algorithm of \cite{Neu1} shows that the signature of the link at infinity is equal to $-5$, so $w(L_\infty)=4$ and $u(L_\infty)=6$.
There is one divorce handle, and indeed $w(L_\infty)=\sum w(\LS_k)-2$.
\end{example}

From Theorem~\ref{mthm2} we can deduce many interesting corollaries. First of all we use it in showing than some
curves with given singularities might not exist. 
The point (a) of the corollary below
is almost a restatement of the result of Petrov \cite{Pet}, which can be interpreted as in \cite{BZ2} as a bound for $k$ with $p=3$.
The point (c) gives the same estimate as in \cite{BZ3}, but we use here only elementary facts, not the
BMY inequality.
\begin{corollary}\label{cormain}
Let $x(t),y(t)$ be polynomials of degree $p$ and $q$ with $p,q$ coprime. Let $C$ be the curve given in parametric form
by 
\begin{equation}\label{eq:Cisrational}
\{w_1=x(t),w_2=y(t),t\in\mathbb{C}\}.
\end{equation}
Assume that the singularity of $C$ at the origin has a branch with singularity $A_{2k}$ $($i.e. $A_{2k}$ is a singularity
of a parametrisation$)$. Then $2k$ is less than or equal to the signature
of the toric knot $T_{p,q}$. In particular
\begin{itemize}
\item[(a)] $k\le q-1-2\intfrac{q}{6}$ if $p=3$;
\item[(b)] $k\le \frac32(q-1)-2\intfrac{q}{4}$ if $p=4$;
\item[(c)] $k\le \sim \frac{pq}{4}$ in general.
\end{itemize}
\end{corollary}
\begin{proof}
Let $L_0$ be the link of singularity of $C$ at $0$. Let $c(L_0)$ be the number of its components. By assumption, one
of its components is a link $T_{2,2k+1}$ with signature $-2k$. By Corollary~\ref{stupid}
\[-\sigma(L_0)\ge 2k+c(L_0)-1.\]
Hence
\[w(L_0)\ge 2k.\]
The link at infinity $L_\infty$ is a knot $T_{p,q}$. Hence $w(L_\infty)=\sigma(L_\infty)=\sigma(T_{p,q})$. 
This, in turn, is computed in Lemma~\ref{torsig}.
The result is then a direct consequence of Theorem~\ref{mthm2}, since $\pg(C)=0$ by assumption (see \eqref{eq:Cisrational}).
\end{proof}
\begin{remark}
Corollary~\ref{cormain}(c) holds even if $p$ and $q$ are not coprime. 
We can compute the signature of the knot at infinity by Lemma~\ref{cable} below.
\end{remark}

The next result is somewhat unexpected, especially if we compare it to \cite[Proposition 87]{Rud} 
stating that no invariant coming from a Seifert matrix of the
knot, including the signature, can tell whether a link is a $\mathbb{C}-$link.
\begin{corollary}\label{ratknot}
If a $\mathbb{C}-$link $L$ with $m$ components bounds an algebraic  curve of geometric genus $\pg$ then
\[-\sigma(L)\ge 2-2m-2\pg.\]
In particular, if a knot bounds a rational curve, its signature is non--positive.
\end{corollary}

Now we can rephrase Theorem~\ref{mthm2} in a Kawauchi--like inequality.
\begin{corollary}\label{betti}
Let $C$ be as in Theorem~\ref{mthm2}. Let $b$ be the first Betti number of $C$ (i.e. the rank of $H_1(C;\mathbb{Q})$. We stress here that we consider
the homology of $C\subset\mathbb{C}^2$, \emph{not} of its compactification in $\mathbb{C}P^2$). Then 
\[|\sigma(L_\infty)-\sum_{k=1}^n\sigma(\LS_k)|\le b+n(L_\infty)-1.\]
\end{corollary}
\begin{proof}
Let $r_k$ be the number of branches of the link $\LS_{k}$ and $d$ be the number of branches at infinity. By Theorem~\ref{mthm2} and
the fact that $w(\LS_k)\ge -\sigma(\LS_k)-(r_k-1)$ we get.
\[-\sigma(L_\infty)-d+n(L_\infty)\ge -\sum\sigma(\LS_k)-\sum(r_k-1)-2(p_g(C)+d-1).\]
Denoting $R=\sum (r_k-1)$ we get
\[\sigma(L_\infty)-\sum\sigma(\LS_k)\le 2p_g+R+d+n(L_\infty)-2=b+n(L_\infty)-1,\]
as $b=2p_g+R+d-1$.
The inequality in the other direction is proved in an identical way, using the invariant
$u$ instead of $w$.
\end{proof}
With not much work, Corollary~\ref{betti} 
can be deduced from \cite{Kaw,Kaw2} (see \cite[Theorem~12.3.1]{Kaw-book}), 
\emph{without ever using the holomorphicity of $C$}. Roughly speaking, we pick a ball $B\subset\mathbb{C}^2$ disjoint from $C$ and pull (by an isotopy)
all the singular points of $C$ inside $B$, so as to get a real surface $C'$ with the property that $C'\cap \partial B$ is a disjoint union of 
links $\LS_1,\dots,\LS_n$. Then $C'\setminus B$ realizes a cobordism between this sum and the link of $C$ at infinity. 
Then \cite[Theorem~12.3.1]{Kaw-book} provides Corollary~\ref{betti}.

The main drawback of that approach is that $C'$ is no longer holomorphic. In short, it works for the signature (and Tristram--Levine
signatures as well),
but if we want at some moment to go beyond and use some more subtle invariant, holomorphicity of $C$ might be crucial. At present we do
not know any such invariant, but we are convinced
that without exploiting thoroughly the holomorphicity of $C$ we cannot get a full understanding of the relation between the link
at infinity and the links of singularities of $C$.

We finish this section by showing a topological proof of Proposition~\ref{nemeth1}. For a convenience of the reader we recall the statement.
\begin{proposition}
Let $f$ be a reduced polynomial in two variables such that the curve $\{f=0\}$ has an isolated
singularity at $(0,0)$. Let $f=f_1\cdot f_2$ be the decomposition of $f$ locally near $(0,0)$, such
that $f_1(0,0)=f_2(0,0)=0$. Let $L$, $L_1$ and $L_2$ be the links of singularities of $\{f=0\}$, $\{f_1=0\}$
and $\{f_2=0\}$ at $(0,0)$ and $\sigma$, $\sigma_1$, $\sigma_2$ its signatures.
Then we have
\[\sigma\le \sigma_1+\sigma_2.\]
\end{proposition}
\begin{proof}
Let $r>0$ be small enough, so that $L=\{f=0\}\cap S(0,r)$ is the link of the singularity of $f$. For a generic vector $v\in\mathbb{C}^2$ sufficiently close to $0$, the intersection
of $S(0,r)$ with $C=C^v=\{F_v=0\}$ is isotopic to $L$, where $F_v(w)=f_1(w)f_2(w-v)$. By definition, $C=C_1\cup C_2$
where 
\[C_1=\{f_1(w)=0\}\cap B(0,r)\text{ and }C_2=\{f_2(w-v)=0\}\cap B(0,r).\]
Let $\varepsilon\ll r$. The link $C\cap S(0,\varepsilon)$ is clearly the link $L_1$ of the singularity given by $\{f_1=0\}$. Consider a change
of the isotopy type of $C\cap S(0,s)$ as $s$ increases from $\varepsilon$ to $r$.

\emph{Claim. There are neither divorce nor fake divorce handles on $C$ for $s\in[\varepsilon,r]$.}

The claim follows from Proposition~\ref{prop:geom}: we put $r_0=\varepsilon$ and $r_1=r$. Then $\pg_1=\pg_0=0$, indeed, the normalization of $C$
is a union of disks. Moreover, in the notation from Proposition~\ref{prop:geom}, $c_1=k_1$ and $c_0=k_0$. In fact, to show $c_0=k_0$ we observe
that $C\cap S(0,\varepsilon)$ is the link of singularity, and both $c_0$ and $k_0$ are the numbers of branches of the singular point. The same argument shows
that $c_1=k_1$ is equal to the number of branches of singularity of $f$ at $(0,0)$. This shows the claim.

Now the Morse theoretical arguments show that
\[w(L)\ge w(L_1)+\sum_{k}w(\LS_k),\]
where we sum over all singular points of $C$, which lie in $B(0,r)\setminus B(0,\varepsilon)$. These singular points are easy to describe. Indeed, there are
no singular points which lie only on $C_1$, there is one
singular point,  at $v$, that lie only on $C_2$ and the corresponding link is the link $L_2$. 
Moreover, there are double points arising as intersections of $C_1$ with $C_2$. The number of these double points can be effectively computed as the
local intersection index of $\{f_1=0\}$ with $\{f_2=0\}$, alternatively as the linking number of $L_1$ with $L_2$,
but we content ourselves by pointing out that for each double point $w(\LS_k)=0$ (see Remark~\ref{rem:blind}). 
Therefore, we get
\[w(L)\ge w(L_1)+w(L_2).\]
And the statement of proposition follows from Lemma~\ref{lem:folklore}, because then $w(L)=-\sigma(L)$, $w(L_1)=-\sigma(L_1)$ and $w(L_2)=-\sigma(L_2)$. 
\end{proof}

\section{Application of Tristram--Levine signatures}\label{secTL}
The notion of signature was generalised by Tristram and Levine \cite{Tr,Le}. The Tristram--Levine signature
turns out to be a very strong tool in the theory of plane algebraic curves. In what follows $\zeta$ will denote
a complex number of module $1$.
\begin{definition}\label{TL}
Let $L$ be a link and $S$ be a Seifert matrix. Consider the Hermitian form
\begin{equation}\label{eq:TLS}
(1-\zeta)V+(1-\bar{\zeta})V^T.
\end{equation}
The \emph{Tristram--Levine} signature $\sigma_\zeta(L)$ is the signature of the above form. The \emph{nullity} $n_\zeta(L)$
is the nullity of above form increased by 1.
\end{definition}
The addition of $1$ is a matter of convention. This makes the nullity additive under disconnected and not connected
sum.
\begin{remark}\label{root}
For a link $L$, let us define $n_0(L)$ as a minimial number such that the $n_0(L)$-th Alexander polynomial is non-zero. 
Let $\Delta_{\min}(L)=\Delta_{n_0(L)}(L)$. Then, 
it is a matter of elementary linear algebra to prove that $n_\zeta(L)\ge n_0(L)+1$ and $n_\zeta(L)>n_0(L)+1$ iff $\Delta_{min}(\zeta)=0$ (we
owe this remark to A.~Stoimenow, see \cite{BN} for a thourough discussion).
\end{remark}
\begin{example}
For $\zeta=-1$ we obtain the classical signature and nullity.
\end{example}
We have, in general, scarce control on the values of $n_\zeta$ if $\zeta$ is a root of the Alexander polynomial. 
However, many interesting
results can be obtained already by studying invariants $\snsz$ and $\nnsz$ when $\zeta$ is not a root of the Alexander polynomial. To
simplify the formulation of these results let us define the functions $\sz$ and $\nz$ as
\begin{equation}\label{eq:sz}
\sz=\begin{cases}
\snsz&\text{ if $\zeta$ is not a root of $\Delta_{min}$}\\
\lim\limits_{\rho\to\zeta^+}\sigma_\rho&\text{ otherwise.}
\end{cases}
\end{equation}
Here $\rho\to\zeta^+$ if we can write $\rho=\exp(2\pi iy)$, $\zeta=\exp(2\pi ix)$ and $y\to x^+$. Similarly we can define $\nz$.
By Remark~\ref{root}, $\nz\equiv n_0(L)+1$, but we keep this function in order to make notation consistent with previous sections.

Tristram--Levine signatures share similar properties as classical signature.
\begin{lemma}[see \cite{Tr,Le}, compare also \cite{St}]\label{tlmur}
Lemma~\ref{mur} holds if we exchange $\sigma(L)$ and $n(L)$ with $\sz(L)$ and $\nz(L)$.
\end{lemma}
Litherland \cite{Li} computes also the signature of torus knot $T_{p,q}$:
\begin{lemma}\label{torsig2}
Let $p,q$ be coprime and $\Sigma$ as in Lemma~\ref{torsig}. Let $\zeta=\exp(2\pi ix)$ with $x\in(0,1)$. Then
\begin{equation}\label{eq:torsig2}
\sz(T_{p,q})=\#\Sigma-2\#\Sigma\cap(x,1+x].
\end{equation}
\end{lemma}
The choice of the closure of the interval $(x,1+x]$ in formula~\eqref{torsig2} agrees with taking the right limit in formula~\eqref{eq:sz}. Indeed,
if $x_k\to x^+$ then the number of points in $\Sigma\cap(x_k,x_k+1]$ converges to the number of points in $\Sigma\cap(x,x+1]$.

The signature of an iterated torus knot can be computed inductively from the result of \cite{Li}.
\begin{lemma}\label{cable}
Let $K$ be a knot and $K_{p,q}$ be the $(p,q)-$cable on $K$. Then for any $\zeta$ we have
\[\sigma_\zeta(K_{p,q})=\sigma_{\zeta^q}(K)+\sigma_\zeta(T_{p,q}).\]
\end{lemma}
This allows recursive computation of signatures of all possible links of unibranched singularities. 
In the general case one uses results of \cite{Neu0,Neu1}.

Because of Lemma~\ref{tlmur} we can repeat the reasoning from Section~\ref{classsig} to obtain a reformulation of Theorem~\ref{mthm2},
Corollary~\ref{ratknot} and Corollary~\ref{betti}.
\begin{theorem}\label{tlbetti}
Let $C$ be an algebraic curve with singular points $z_1,\dots,z_n$, with links of singularities $\LS_1,\dots,\LS_n$. 
Let $L_\infty$
be the link of $C$ at infinity. Let also $b$ be the first Betti number of $C$. Then
\begin{equation}\label{eq:tlbetti}
\left|\sz(L_\infty)-\sum\sz(\LS_k)\right|\le b+n_0(L_\infty).
\end{equation}
\end{theorem}
The proof goes along the same line as the proof of Corollary~\ref{betti}. We introduce the quantities $w_\zeta=-\sz(L)+\nz(L)-c(L)$
and $u_\zeta=-\sz(L)-\nz(L)+c(L)$ and study their changes on crossing different singular handles.
We remark only that $\nz(L_\infty)=n_0(L_\infty)+1$.

Using the same argument as in Proposition~\ref{prop:geom} we obtain a result which relates the signatures at two intermediate steps.

\begin{proposition}\label{aftertlbetti}
For any generic parameter $\xi$, let 
$r_0$ and $r_1$ be two non-critical parameters. For $i=0,1$ let $L_i$, $c_i$  be, respectively, the link $C\cap S(\xi,r_i)$ and its
number of components. Let $\Delta\pg$ be the difference of genera of $C\cap B(\xi,r_1)$ and $C\cap B(\xi,r_0)$ and $\Delta k$ the difference
between number of connected components of corresponding normalizations. We have then
\begin{align*}
w_\zeta(L_1)-\sum w_\zeta(\LS_k)-w_\zeta(L_0)&\ge -2(\Delta\pg+c_1-c_0-\Delta k),\\
-(u_\zeta(L_1)-\sum u_\zeta(\LS_k)-u_\zeta(L_0))&\ge -2(\Delta\pg+c_1-c_0-\Delta k),
\end{align*}
where we sum only over those critical points that lie in $B(\xi,r_1)\setminus B(\xi,r_0)$.
\end{proposition}

Corollary~\ref{ratknot} generalises immediately to the following, apparently new result.
\begin{lemma}\label{ratknot2}
If $K$ is a $\mathbb{C}-$knot bounding a rational curve, then $\sz(K)\le 0$ for any $\zeta$.
\end{lemma}

Another application of Theorem~\ref{tlbetti} is in the classical problem of bounding the number of cusps of a plane curve of degree $d$,
see \cite{Hir} for the discussion of this problem. Our result is a topological proof of Varchenko's bound.
\begin{corollary}\label{cor:2372}
Let $s(d)$ be a maximal number of $A_2$ singularities on an algebraic curve in $\mathbb{C}P^2$ of degree $d$. Then
\[\limsup \frac{s(d)}{d^2}\le \frac{23}{72}.\]
\end{corollary}
\begin{proof}[Proof (sketch)]
Let $C$ be a curve of degree $d$ in $\mathbb{C}P^2$. Let us pick up a line $H$ intersecting $C$ in $d$ distinct points. We chose
an affine coordinate system on $\mathbb{C}P^2$ such that $H$ is the line at infinity. Let $C_0$ be the affine part of $C$. Then $C_0$ can be defined as a zero
set of a polynomial $F$ of degree $d$.
Let $z_1,\dots,z_s$ be the singular points of $C_0$ of type $A_2$.

\underline{\emph{Case 1.}} $C_0$ has no other singular points.

Then, by the genus formula, $b_1(C_0)=d^2-2s+O(d)$. Let us take $\zeta=e^{\pi i/6}$. Then $\sz(\LS_i)=2$. On the other hand, the link of $C_0$ at infinity
is toric link $T_{d,d}$ and its signature 
\[\sz(T_{d,d})=2d^2\cdot\frac16\left(1-\frac16\right)+O(d)=\frac{5}{18}d^2+O(d).\]
(For $\zeta=e^{2\pi ix}$ we have asymptotics $\sz(T_{d,d})=2d^2x(1-x)+O(d)$ by results \cite{Neu0,Neu1}.)
Then \eqref{eq:tlbetti} provides
\[2s-\frac{5}{18}d^2\le d^2-2s+O(d).\]

\smallskip
\underline{\emph{Case 2.}} $C_0$ has other singular points.
Let $\xi\in \mathbb{C}^2$ be a generic point of $\mathbb{C}^2$.
and let $r_\infty$ be sufficiently large, so that the intersection
of $C_0$ with a sphere $S(\xi,r_\infty)$ is transverse. Let $G$ be a generic polynomial of very high degree vanishing
at each of $z_k$ with up to order at least $4$ (i.e. generic among polynomials sharing this property). For $\varepsilon>0$ small enough
this guarantees that the curve 
\[C_\varepsilon=\{F+\varepsilon G=0\}\]
has singularities of type $A_2$ at each $z_k$, is smooth in $B(\xi,r_\infty)$ away from $z_k$'s and
its intersection with the sphere $S(\xi,r_\infty)$ is the same as the intersection
of $C_0$. Now we can repeat the proof in Case~1.
\end{proof}
The above estimate is very close to the best known to the author, that the limit is bounded from above by $(125+\sqrt{73})/432$
(\cite{Lan}). 

Theorem~\ref{tlbetti} can be used together with results (especially Lemma~3 and Theorem~3) in \cite{Li}. 
We can get another proof of classical Zajdenberg--Lin theorem (see \cite{LZ}), if we put $b=0$ (we defer the details to a subsequent paper). 
It is, presumably, possible to go beyond this theorem and classify all plane curves with small first Betti number
(compare \cite{BZ0} and \cite{BZ1}). We can
also hope to prove some results concerning the maximal possible number of singular point of the algebraic curve with given
first Betti number, the problem that is known as the Lin conjecture.

\begin{acknowledgements}
The author is very grateful to A.~N\'emethi for various discussions on the subject and providing a proof of Proposition~\ref{nemeth1}. He
also would like to express his thanks to H.~\.Zo\l{}\c{a}dek for carefully reading the manuscript at early stage of
its preparation.
The author would also like to thank A. Stoimenow and P. Traczyk for patiently explaining some elements of knot
theory obscure to the author. He is also grateful to A. P\l{}oski for many stimulating discussions and to L.~Rudolph for his interest in this work.
\end{acknowledgements}

\end{document}